\global\newcount\secno \global\secno=1
\global\newcount\propno \global\propno=1 
\global\newcount\eqnum \global\eqnum=1 

\def\prop#1{\xdef #1{\the\secno.\the\propno} 
            \global\advance\propno by1  } 
\def\Eqno#1{\xdef #1{\the\secno.\the\eqnum} \eqno(\the\secno.\the\eqnum)
            \global\advance\eqnum by1  
            \eqlabeL #1}


\def\proplabeL#1{} 
\def\eqlabeL#1{} 

\def\mapright#1{\smash{\mathop{\longrightarrow}\limits^{#1}}}

\def\mapdown#1{\Big\downarrow\rlap{$\vcenter{\hbox{$\scriptstyle#1$}}$}}
\def\mapdownleft#1{\Big\downarrow\llap{$\vcenter{\hbox{$\scriptstyle#1$}}\;\;$}}

\def\simrightarrow{\smash{\mathop{\rightarrow}\limits^{\sim}}} 
\def\vM{{\mathaccent20 M}}
\def\vtau{{\mathaccent20 \tau}}
\def\vX{{\mathaccent20 X}}
\def\ve{{\mathaccent20 e}}
\def\vf{{\mathaccent20 f}}
\def\vB{{\mathaccent20 B}}
\def\vK{{\mathaccent20 K}}

\def\preMainI{1.6}
\def\preMainII{1.9}
\def\preMainIII{1.17}

\overfullrule=0pt

\input epsf
\input amstex.tex 
\def\Proclaim#1#2{ \prop{#2} \proclaim{#1 #2} \proplabeL{#2} }
\documentstyle{amsppt} 
\topmatter
\title
Autoequivalences of Derived Category of A K3 surface 
and Monodromy Transformations 
\endtitle
\author
Shinobu Hosono, Bong H. Lian, Keiji Oguiso, Shing-Tung Yau
\endauthor
\address 
\endaddress

\abstract 
We consider autoequivalences of the bounded derived category of 
coherent sheaves on a K3 surface. 
We prove that the image of the autoequivalences has index at most two 
in the group of the Hodge isometries of the Mukai lattice. 
Motivated by homological mirror symmetry we also consider the mirror 
counterpart, i.e. symplectic version of it. In the case of $\rho(X)=1$, 
we find an explicit formula which reproduces the number of Fourier-Mukai 
(FM) partners from the monodromy problem of the mirror K3 family. 
We present an explicit example in which a monodromy action does not 
come from an autoequivalence of the mirror side. 
\endabstract

\leftheadtext{S. Hosono, B.H. Lian, K. Oguiso and S.-T. Yau}
\rightheadtext{Autoequivalence of Derived Category of K3 surfaces}
\endtopmatter

\head 
{ Table of Contents }
\endhead
{  
\leftskip1cm\rightskip1cm
\item{\S 0} Introduction  -- Motivation and Backgrounds 
\item{\S 1} Statements of main results
\item\item{(1-1)} Autoequivalences and 
 $\text{Im}(\text{Auteq} D(X) \rightarrow O_{Hodge}(\tilde H(X,\bold Z)) )$
\item\item{(1-2)} Symplectic diffeomorphisms and surjectivity to $O^+(T(Y))^*$
\item\item{(1-3)} Mirror symmetry of marked $M$-polarized K3 surfaces 
\item\item{(1-4)} Fourier-Mukai partners and monodromy of the mirror family 
\item{\S 2} Autoequivalences and Proof of Theorem {\preMainI}
\item\item{(2-1)} Various autoequivalences 
\item\item{(2-2)} FM transforms on a K3 surface
\item\item{(2-3)} Proofs of Theorem {\preMainI} 
\item{\S 3} Symplectic mapping class group and Proof of Theorem {\preMainII}
\item{\S 4} Monodromy group, FM partners, and Proof of Theorem {\preMainIII}
\item{\S 5} Mirror family of a K3 surface with $\text{deg}(X)=12$
\par}

\document

\vskip1.5cm

\head
{\S 0. Introduction -- Motivation and Backgrounds}
\endhead 

Our main results are Theorems {\preMainI}, {\preMainII}, {\preMainIII} 
and Proposition 5.8. 

\vskip0.3cm

In this paper we study the bounded derived category of K3 surfaces 
motivated by homological mirror symmetry of K3 surfaces.

Homological mirror symmetry due to Kontsevich [Ko] is based on 
homological and algebraic aspects of manifolds where one considers 
certain derived categories of manifolds, i.e. 
the bounded derived category $D(X)$ of a projective variety $X$, 
and the bounded derived category $D Fuk(\vX,\beta)$ of the Fukaya category 
$Fuk(\vX,\beta)$ for the mirror $\vX$ with its symplectic structure $\beta$. 
(See [FO3], [Fu] and reference therein.)  
The homological mirror symmetry conjecture says 
that {\it when $X$ and $\vX$ form a mirror pair, there should be 
a (not necessarily canonical) exact equivalence} 
$$
D(X) \cong D Fuk(\vX,\beta),
$$
where $\beta$ is a generic symplectic structure and a generic complex 
structure is assumed for the left hand side. 
Despite its homological and algebraic nature, 
homological mirror symmetry bears a close 
relationship with geometric mirror symmetry due to [SYZ] which is based  
on the symplectic geometry of the underlying $C^\infty$-manifolds. 
See also [Mo2] and reference therein for details. 

The following mathematical problems are important for understanding 
fully the homological mirror symmetry:  

\item{a)} determine the group of autoequivalences $D(X) \rightarrow D(X)$,
\item{b)} determine the set of varieties $Y$ s.t. $D(Y) \cong D(X)$, 
i.e. FM partners of $X$.

\noindent
Recent results by Orlov[Or2] solve both problems when $X$ is an abelian 
variety. For K3 surfaces, the problem b) has been studied  
in detail [Or2],[BM],[Og2].

In this paper, we first consider 
the natural map from the autoequivalences 
to the Hodge isometries of the Mukai lattice and prove that the image 
is a subgroup of index at most two ({\bf Theorem {\preMainI}}). 
Next we consider the mirror counterpart of the 
problem a), i.e. the group of autoequivalences of 
$D Fuk(\vX,\beta)$. However at this moment not much is known 
about this side.  So we consider 
the group $\text{Symp}(\vX,\beta)$ of the (cohomological) symplectic 
diffeomorphisms, more precisely the (cohomological) symplectic mapping 
class group $\pi_0 Symp(\vX,\beta)=Symp(\vX,\beta)/Symp^0(\vX,\beta)$,  
as a natural substitute for $\text{Auteq} D Fuk(\vX,\beta)$ (Definition 1.7). 
Under this substitution we prove that the pullback of the differential forms 
induces a surjective map from $\text{Symp}(\vX,\beta)$ to an orthogonal 
group $O^+(T(\vX))^*$ of the transcendental lattice 
({\bf Theorem {\preMainII}}).

Second, we consider a mirror family in 
the sense of Dolgachev [Do] and define 
a monodromy group $\Cal M(\vX)$ of the family. Then using the surjective 
map in Theorem {\preMainII} we define a monodromy representation of the group 
$\pi_0 \text{Symp}(\vX,\beta)$ and consider its image 
$\Cal{MS}(\vX)=O^+(T(\vX))^*$ in $\Cal{M}(\vX)$. 
Surprisingly, for K3 surfaces with 
$\rho(X)=1$, we find that the group index 
$[ \Cal{M}(\vX): \Cal{MS}(\vX)]$ coincides with the number of 
FM (Fourier-Mukai) partners obtained in [Og2] ({\bf Theorem {\preMainIII}}). 

Finally, we present explicit calculations of the group index of the 
monodromy representation in the first non-trivial case $X$ of $\text{deg}(X)
=12$. 
That the monodromy of the mirror family are realized by  
autoequivalences in mirror symmetry was first suggested by Kontsevich 
(see e.g. [Mo2]) and this has been studied further for example in 
[Hor],[ST],[ACHY]. Our example shows that this need not be the case 
when we have non-trivial FM partner ({\bf Proposition 5.8}). 
Our example has appeared in [LY1][LY2], and also studied in detail 
in completely different context [PS][BP] (before the mirror symmetry).  

\vskip0.5cm
After we have posted the preliminary version of this paper, 
B. Szendr\"oi kindly pointed out to us a mistake in the 
preliminary version. We are also grateful to him for informing us of his 
work in mirror symmetry of K3 surfaces. (See section 5 of [Sz].) 

We also thank D. Huybrechts for informing us of his student, 
D. Ploog's work on the map $\text{Auteq} D(X) \rightarrow 
O_{Hodge}(\tilde H(X,\bold Z))$. 

\vfill
\eject

\head
{Acknowledgement}
\endhead 
The first and the third named authors  
would like to thank the Department of Mathematics of 
Harvard University for the hospitality during their stay. They would like 
to thank the Education Ministry of Japan and the Harvard University 
for their financial support during their stay. The second named author 
would like to thank D. Ruberman for helpful discussions. He is supported by 
NSF grant DMS-0072158.  
The authors would like to thank K. Fukaya for informing us 
of the latest reference about $D Fuk(\vX,\beta)$.

\head
{\S 1. Statements of main results}
\endhead 

\global\secno=1  
\global\propno=1 
\global\eqnum=1 

\noindent
(1-1) Autoequivalences and 
 $\text{Im}(\text{Auteq} D(X) \rightarrow O_{Hodge}(\tilde H(X,\bold Z)) )$
\vskip0.2cm

Let $X$ be a K3 surface, a smooth projective surface over $\bold C$ with 
$\Cal O_X(K_X)\cong \Cal O_X$ and $h^1(\Cal O_X)=0$. We denote by $(*,**)$ the 
symmetric bilinear form on $H^2(X,\bold Z)$ given by the cup product. Then 
$(H^2(X,\bold Z), (*,**))$ is an even unimodular lattice of signature 
$(3,19)$. This lattice is isomorphic to the K3 lattice $\Lambda_{\text{K3}}:=
E_8(-1)^{\oplus 2}\oplus U^{\oplus 3}$ where $U$ is the hyperbolic lattice 
( an even unimodular lattice of signature (1,1)). 
An isomorphism $\tau: H^2(X,\bold Z)
\rightarrow \Lambda_{\text{K3}}$ is called a marking and a pair 
$(X,\tau)$ is called a marked K3 surface. We denote by 
$NS(X)=\text{Pic}(X)$ the N\'eron-Severi lattice of $X$ and by 
$\rho(X)$ the Picard number, i.e. the rank of $NS(X)$. The lattice $NS(X)$ 
is primitive in $H^2(X,\bold Z)$ and has signature $(1,\rho(X)-1)$. We call 
the orthogonal lattice $T(X):=NS(X)^\perp$ in $H^2(X,\bold Z)$ the 
transcendental lattice. $T(X)$ is primitive 
in $H^2(X,\bold Z)$ and has signature $(2,20-\rho(X))$. We denote by 
$\omega_X$ a nowhere vanishing holomorphic two form. Then the natural 
inclusion 
$$
\bold C \omega_X \oplus \bold C \bar \omega_X \subset T(X)\otimes \bold C
$$
defines a Hodge structure of weight 2 on $T(X)$.

Based on the original work by Mukai[Mu2], Orlov[Or1] showed that Fourier-Mukai 
transform on the bounded derived category of coherent sheaves induces 
a Hodge isometry of the Mukai lattice. Let us summarize basic definitions.

\Proclaim{Definition}{\MukaiLattice} 
For a K3 surface $X$, we define the {\it Mukai lattice} 
to be a lattice 
$$
\tilde H(X,\bold Z):= \big( H^0(X,\bold Z) \oplus H^2(X,\bold Z) \oplus 
H^4(X,\bold Z), \langle\;,\; \rangle \big) \;\;,
$$
with its non-degenerate bilinear form 
$\langle (a,b,c),(a',b',c')\rangle :=-(a,c')-(c,a')+(b,b')$ using  
the cup product $(\;,\;)$ of the cohomology ring $\oplus H^i(X,\bold Z)
=H^0(X,\bold Z)\oplus H^2(X,\bold Z) \oplus H^4(X,\bold Z)$. Then 
there is an isomorphism as lattice $\tilde H(X,\bold Z) 
\cong U \oplus \Lambda_{K3}$, where $U$ is the hyperbolic lattice 
and $\Lambda_{K3}=U^{\oplus 3} \oplus E_8(-1)^{\oplus 2}$.
\endproclaim

\Proclaim{Definition}{\HodgeIso}  An isometry 
of the Mukai lattices $\varphi: \tilde H(X,\bold Z)\rightarrow \tilde 
H(Y,\bold Z)$ is called Hodge isometry if it satisfies 
$
\varphi(\bold C \omega_X) = \bold C \omega_Y, 
$ 
and is denoted by 
$$
\varphi:(\tilde H(X,\bold Z),\bold C \omega_X) \rightarrow 
(\tilde H(Y,\bold Z),\bold C \omega_Y).
$$ 
For $Y=X$, we define the group of 
Hodge isometries 
$$
O_{Hodge}(\tilde H(X,\bold Z)):=
\{ \varphi: (\tilde H(X,\bold Z),\bold C \omega_X) \rightarrow 
            (\tilde H(X,\bold Z),\bold C \omega_X) \}.
$$
\endproclaim

\vskip0.2cm

For a smooth projective variety $X$ we denote by $D(X):=D_{coh}^b(X)$ 
the bounded derived category of coherent sheaves on $X$, which hereafter 
will be called simply the derived category of $X$. $ObD(X)$ then consists of 
the bounded complexes of coherent sheaves on $X$. 
$D(X)$ is naturally regarded as a triangulated category. 
(See [GM] for details.)  In this paper, 
by a functor $F: D(X) \rightarrow D(Y)$ we always mean a functor 
as a triangulated category, i.e. a functor which commutes with 
the shift functor and preserves the distinguished triangles.

A functor $F: D(X) \rightarrow D(Y)$ is called an {\it equivalence} if 
there exists a functor $G:D(Y)\rightarrow D(X)$ for which we have the 
relations $G\circ F \cong \text{id}_{D(X)}$ and 
$F\circ G \cong \text{id}_{D(Y)}$ as functors. $G$ is called {\it 
quasi-inverse } of $F$. 
The isomorphism class of $G$ (as a functor) is uniquely determined by the 
isomorphism class of $F$, although $G$ itself is not uniquely determined. 
The group of the isomorphism classes 
of (self-)equivalences $F:D(X) \rightarrow D(X)$ is called an 
{\it autoequivalence} of $D(X)$ and will be denoted by $\text{Auteq} D(X)$. 
Note that $\text{Auteq} D(X)$ is a set and hence a group because of the 
finiteness condition on the coherent sheaves. 

For smooth projective varieties $X,Y$ and $\Cal E \in D(X \times Y)$, 
we consider a functor
$\Phi^{\Cal E}_{X\rightarrow Y}: D(X) \rightarrow D(Y)$ defined by
$$
\Phi^{\Cal E}_{X\rightarrow Y}(\Cal X)=
\bold R \pi_{Y*}( \Cal E \otimes^{\hskip-7pt \,^\bold L} \bold L \pi_{X}^* 
\Cal X )
\Eqno\FM
$$ 
where $\pi_X:X\times Y \rightarrow X$ and $\pi_Y:X\times Y \rightarrow Y$ 
are the natural projections; $\pi_X(a,b)=a, \pi_Y(a,b)=b$. 
In this paper we abbreviate the functor (\FM) as  
$\Phi^{\Cal E}_{X\rightarrow Y}(\Cal X)=\pi_{Y*}(\Cal E \otimes \pi^*_{X} 
\Cal X)$ and also write $\Phi^{\Cal E}$ if $X=Y$. 

In what follows we will be mainly concerned with the case where 
the functor $\Phi^{\Cal E}_{X\rightarrow Y}$ gives an equivalence. In 
this case, this functor is called a {\it Fourier-Mukai(FM) transform} and 
$\Cal E \in D(X\times Y)$ is called its {\it kernel}. The following 
fundamental result by Orlov allows us to represent an equivalence in the 
form of FM transform:

\Proclaim{Theorem}{\Or} {\bf ([Or1])} 
For smooth projective varieties $X,Y$ and 
an equivalence $\Phi_{X\rightarrow Y}: D(X) \rightarrow D(Y)$, 
there exists $\Cal E \in D(X\times Y)$ such that 
$\Phi_{X\rightarrow Y}=\Phi^{\Cal E}_{X\rightarrow Y}$. Moreover 
$\Cal E$ is determined by $\Phi$ unique up to isomorphism as an object in 
$D(X \times Y)$.  \qed
\endproclaim

For an element $\Cal F =(\Cal F^\bullet) \in ObD(X)$ 
we define the chern character $ch(\Cal F)= \sum_i (-1)^i$ $ch(\Cal F^i) 
\in \tilde H(X,\bold Q):=\oplus_k H^{2k}(X,\bold Q)$. 
This is well-defined by the definition of $D(X)$. Furthermore by 
the uniqueness of the kernel 
$\Cal E \in D(X \times Y)$ up to isomorphism in Theorem {\Or}, the chern 
character $ch(\Cal E)$ does not depend 
on the choice of the kernel representing an equivalence $\Phi$. 

Now let us restrict our attention to K3 surfaces. The following is a 
fundamental Theorem essentially due to Mukai:

\Proclaim{Theorem}{\D(X)} Let $X$ be a K3 surface and $Y, Z$ 
smooth projective manifolds. 
\item{1)} If $\Phi: D(X) \rightarrow D(Y)$ is an equivalence, then $Y$ is a 
K3 surface. 
\item{2)} For each $\Cal X \in D(X)$, $ch(\Cal X) \in \tilde H(X,\bold Z)$.
\item{3)} An equivalence $\Phi^{\Cal E}_{X\rightarrow Y}:D(X) 
\rightarrow D(Y)$ induces a Hodge isometry 
$$
f^{\Cal E}_{X\rightarrow Y}: \tilde H(X,\bold Z) \rightarrow 
\tilde H(Y,\bold Z)
$$
defined by $f^{\Cal E}_{X\rightarrow Y}(x)=
\pi_{Y*}(\Cal Z \cdot \pi_X^* x)$, where 
$$
\Cal Z:=\pi_X^*(\sqrt{td_X})ch(\Cal E) \pi_Y^*(\sqrt{td_Y}) \in \tilde 
H(X\times Y, \bold Z) \;\;.
$$
\item{4)} If both $\Phi^{\Cal E}_{X\rightarrow Y}: D(X)\rightarrow D(Y)$ 
and $\Phi^{\Cal E'}_{Y\rightarrow Z}: D(Y)\rightarrow D(Z)$ are equivalence, 
then $f^{\Cal E''}_{X\rightarrow Z}=f^{\Cal E'}_{Y\rightarrow Z} \circ 
f^{\Cal E}_{X\rightarrow Y}$ holds as an isometry $\tilde H(X,\bold Z) 
\rightarrow \tilde H(Z,\bold Z)$, where $\Cal E'' \in D(X \times Z)$ 
is a kernel representing the equivalence $\Phi^{\Cal E'}_{Y\rightarrow Z}
\circ \Phi^{\Cal E}_{X\rightarrow Y}: D(X) \rightarrow D(Z)$. 
\endproclaim

Property 1) is well-known (see for example [Mu1],[BM]). 
Property 2) and the fact that ${\Cal Z} \in \tilde H(X \times Y,\bold Z)$ 
in 3) are due to the even property of $H^2(X,\bold Z)$ for K3 surfaces. 
Property 3) follows from the following facts: 
the result by Mukai (Theorem 4.9,[Mu2]) shows that 
$f^{\Cal E}_{X\rightarrow Y}$ is an isomorphism, 
and further it preserves the bilinear form 
$\langle\;,\;\rangle$ (Lemma 4.7,[Mu2]). Since the kernel 
$\Cal E \in D(X\times Y)$ is algebraic, this map 
preserves the Hodge decomposition and therefore $\bold C \omega_X$ 
is mapped to $\bold C \omega_Y$. The property 4) follows from 
the projection formula and the Grothendieck-Riemann-Roch theorem. (See an 
argument [Mu2, pp. 382--383].)

\Proclaim{Corollary}{\cor} Let $X$ be a K3 surface. An autoequivalence 
$\Phi^{\Cal E}:=\Phi^{\Cal E}_{X \rightarrow X}: D(X) \rightarrow D(X)$ 
induces a Hodge isometry $f^{\Cal E}:=f^{\Cal E}_{X\rightarrow X}: 
\tilde H(X,\bold Z)\rightarrow \tilde H(X,\bold Z)$ which makes the 
following diagram commutative: 
$$
\matrix D(X) & \mapright{\Phi^{\Cal E}} & D(X) \cr
    \mapdownleft{ch(\cdot)\sqrt{td_X}} & & \mapdown{ch(\cdot)\sqrt{td_X}} \cr
      \tilde H(X,\bold Z) & \mapright{f^{\Cal E}} & \tilde H(X,\bold Z) \cr
\endmatrix
$$
We write $f^{\Cal E}=ch(\Phi^{\Cal E})$. Then we have 
$ch(\Phi^{\Cal E_1}\circ \Phi^{\Cal E_2})=
ch(\Phi^{\Cal E_1})\circ ch(\Phi^{\Cal E_2})$, i.e.  we have a group 
homomorphism:
$$
ch: \text{Auteq} D(X) \rightarrow O_{Hodge}(\tilde H(X,\bold Z))
\Eqno\AutoToHodge
$$
\endproclaim

\Proclaim{Theorem}{\MainThI} {\bf (Main Theorem 1)} 
Let $\iota_2$ be 
the involution which maps $(a,b,c) \in \tilde H(X, \bold Z)$ to 
$(a, -b, c)$. Then    
$$
\langle ch(Auteq D(X)), \iota_2 \rangle = O_{Hodge}(\tilde H(X,\bold Z)).  
$$
In particular, the subgroup $ch(\text{Auteq} D(X))$ has index at most two in 
$O_{Hodge}(\tilde H(X,\bold Z))$.
\endproclaim

Our proof of Theorem {\MainThI} is an easy combination of the arguments in 
[Mu2],[Or],[BM] together with the reflection functors by  
spherical objects defined in [ST]. We present its full details in the 
next section.

\vskip0.5cm
\noindent
(1-2) Symplectic diffeomorphisms and surjectivity to $O^+(T(Y))^*$
\vskip0.2cm

Let $Y$ be a  K3 surface and $A(Y) \subset NS(Y)\otimes \bold R$ 
the ample cone, i.e.  an open convex cone in $NS(Y)\otimes \bold R$ generated 
by the ample classes. We note that due to [Ya] and the real Nakai-Moishezon 
criterion [CP](, see also [To]), any class in $A(Y)$ is represented by 
a unique Ricci-flat K\"ahler form on $Y$.  

\Proclaim{Definition}{\sympY} Let $\kappa_Y \in A(Y)$, and $Diff(Y)$ be 
the diffeomorphism group of the underlying real $C^\infty$-four manifold $Y$. 
\item{1)} The pair $(Y,\kappa_Y)$ or simply $\kappa_Y$ is called a 
(cohomological) symplectic structure  of the underlying real 
$C^\infty$-four manifold of $Y$. 
\item{2)} The group $Symp(Y,\kappa_Y):=\{ g \in Diff(Y) \; | \; 
g^* \kappa_Y=\kappa_Y \}$ is called the (cohomological) symplectic 
diffeomorphism group 
of $(Y,\kappa_Y)$, where $g^*:H^2(Y,\bold Z) \rightarrow H^2(Y,\bold Z)$ 
is the isometry induced by $g$. 
\item{3)} The group $\pi_0 Symp(Y,\kappa_Y):=Symp(Y,\kappa_Y)/
Symp^0(Y,\kappa_Y)$ is called the symplectic mapping class group of 
$(Y,\kappa_Y)$, where $Symp^0(Y,\kappa_Y)$ is the connected component of 
the identity. 
\item{4)} A symplectic structure $(Y,\kappa_Y)$ is said to be {\it generic} 
if every primitive lattice $P \subset H^2(Y,\bold Z)$ with 
$ \kappa_Y \in P\otimes \bold R$  contains $NS(Y)$.  
In other words, $(Y,\kappa_Y)$ is generic if it characterizes 
$NS(Y)$ as the minimal primitive sublattice in $H^2(Y,\bold Z)$ which 
contains $\kappa_Y$ after the tensor product $\otimes \bold R$. 
\qed
\endproclaim

\noindent
{\bf Remark.} 1) If we choose an integral basis $e_1,\cdots,e_\rho$ of $NS(Y)$ 
and represent $\kappa_Y=\sum \alpha_i e_i \;(\alpha_i \in \bold R)$, 
then the condition of $\kappa_Y$ generic is satisfied if and only if 
$\alpha_1, \cdots, \alpha_\rho$ are linearly independent over $\bold Q$. 
Therefore the set of generic symplectic structures $\kappa_Y$ is the 
complement of countably many proper hyperplanes in $A(Y)$.  

\noindent
2) Our definition of $Symp(Y,\kappa_Y)$ is {\it cohomological} in nature. 
It need not preserve a symplectic form representing the class $\kappa_Y$. 
Since the mirror map relates complex structures with positive classes 
(cf. [Do]), this group seems more natural than the group of symplectomorphisms 
at least in our context.  \qed 

\vskip0.5cm

Let $L$ be a lattice of signature $(p,q)$ and consider its isometry 
group $O(L)$. We denote by $O^+(L)$ the subgroup consisting those 
isometries which preserve the orientation of the maximal positive definite 
subspaces in $L\otimes \bold R$. 
More precisely, $g \in O^+(L)$ if and only if $\pi \circ g$ induces an 
orientation preserving isomorphism for every positive $p$-space $H$ in 
$L\otimes \bold R$, where $\pi: L\otimes \bold R \rightarrow H$ is the 
orthogonal projection. 
For the transcendental lattice $T(Y)$, we consider the group of 
isometries $O(T(Y))$ and the natural orthogonal representation of $O(T(Y))$ 
on its discriminant group $A_{T(Y)}:=T(Y)^*/T(Y)$ with its 
($\bold Q/\bold Z$-valued) bilinear form 
naturally induced from that of $T(Y)$. Then we define  
$$
O(T(Y))^*:=\text{Ker}( O(T(Y)) \rightarrow O(A_{T(Y)}) ) 
\Eqno{\Ostar}
$$
and set $O^+(T(Y))^* := O(T(Y))^* \cap O^+(T(Y))$. 

\Proclaim{Theorem}{\Donaldson} {\bf ([Don, Section VI])} For a K3 surface $Y$, 
we have: 
\item{1)} If $g \in Diff(Y)$, then 
$g^*|_{H^2(Y,\bold Z)} \in O^+(H^2(Y,\bold Z))$.
\item{2)} The natural map $Diff(Y) \rightarrow O^+(H^2(Y,\bold Z))$ 
given by $1)$ is surjective. 
\qed
\endproclaim

Using this, we prove in section 3 the following symplectic analogue of 
Theorem {\MainThI}: 

\Proclaim{Theorem}{\MainThII} {\bf (Main Theorem 2)} 
Let $(Y,\kappa_Y)$ be a generic symplectic structure 
on $Y$. Then:
\item{1)} If $g \in Symp(Y,\kappa_Y)$, then $g^*(NS(Y))=NS(Y)$, $g^*(T(Y))
=T(Y)$ and $g^*|_{NS(Y)}$ $=\text{id}_{NS(Y)}$, $g^*|_{T(Y)} \in O^+(T(Y))^*$. 
\item{2)} The natural map $\pi_0 Symp(Y,\kappa_Y) 
\rightarrow O^+(T(Y))^*$ induced by 1) is surjective. 
\endproclaim

\noindent
{\bf Remark.} The kernel $O(T(Y))^*$ is a subgroup of $O(T(Y))$ consisting 
those elements $\phi_{T(Y)}$ such that 
$(\phi_{T(Y)},\text{id}_{\widetilde{NS}(Y)}) \in O(T(Y)) \times 
O(\widetilde{NS}(Y))$ extend to isometries in $O(\tilde H(Y,\bold Z))$. 
In the context of the mirror symmetry, we will identify the K3 surface 
$Y$ with the mirror $\vX$ of a K3 surface $X$. There the group $O^+(T(Y))^*$ 
will be identified with a subgroup of $O_{Hodge}(\tilde H(X,\bold Z))$.
\qed

\vskip0.5cm
\noindent
(1-3) Mirror symmetry of marked $M$-polarized K3 surfaces 
\vskip0.2cm

Our results (Theorem {\MainThI} and Theorem {\MainThII}) have clear 
interpretations in terms of homological mirror symmetry of K3 surfaces. 
In order to describe them, we discuss mirror symmetry of   
marked $M$-polarized K3 surfaces following Dolgachev[Do]. 

Let us consider a lattice $M$ of signature $(1,t)$ and assume a primitive 
embedding $\iota_M: M \hookrightarrow \Lambda_{\text{K3}}$. We fix this 
embedding $\iota_M$ and identify $M$ and $\iota_M(M)$ in 
$\Lambda_{\text{K3}}$. Then a pair $(X,\tau)$ of a K3 surface $X$ and a 
marking $\tau: H^2(X,\bold Z)\rightarrow \Lambda_{\text{K3}}$ is called a 
{\it marked $M$-polarized K3 surface} if $\tau^{-1}(M) \subset NS(X)$. 
We call a K3 surface X 
$M$-{\it polarizable} if there is a marking $\tau: H^2(X,\bold Z) \rightarrow 
\Lambda_{\text{K3}}$ such that $(X,\tau)$ is a marked $M$-polarized K3 
surface. Two marked $M$-polarized K3 surfaces $(X,\tau)$ and 
$(X',\tau')$ are said isomorphic if there exists an isomorphism $\varphi:X 
\rightarrow X'$ such that $\tau'=\tau\circ \varphi^*$. 

Let $(X,\tau)$ be a marked $M$-polarized K3 surface and $\omega_X$ 
be a nowhere vanishing holomorphic two form. Since 
$NS(X)=H^{1,1}(X) \cap H^2(X,\bold Z)=(\bold C \omega_X)^\perp \cap 
H^2(X,\bold Z)$, the line $\tau(\bold C \omega_X)$ is always 
orthogonal to $M\otimes \bold C$. Also since $(\omega_X,\omega_X)
=0$, $(\omega_X,\bar \omega_X) >0$, the line $\tau(\bold C \omega_X)$ lies 
in the period domain 
$$
\Omega(M^\perp):= \{ \bold C \omega \in \bold P(M^\perp \otimes \bold C) 
\; \vert \; (\omega,\omega)=0, (\omega,\bar \omega)>0 \;\}
$$
for the orthogonal lattice $M^\perp$ in $\Lambda_{\text{K3}}$. The point 
$\tau(\bold C \omega_X)$ is called 
{\it period point} of $(X,\tau)$. A local family of marked 
$M$-polarized K3 surfaces is a family $f:{\frak X} \rightarrow \Cal B$ of K3 
surfaces together with a trivialization (i.e. a marking) $\tau=\cup_{t \in 
\Cal B}\tau_t:R^2f_* \bold Z_{\frak X} \simrightarrow 
\Lambda_{\text{K3}}\times {\Cal B}$ such that $\tau^{-1}(M \times \{ t\}) 
\subset NS({\frak X}_t)$. Then the period map, 
which sends each ${\frak X}_t (t \in {\Cal B})$ to its 
period points in $\Omega(M^\perp)$, defines a holomorphic map  
from $\Cal B$ to $\Omega(M^\perp)$. Due to the surjectivity of the period 
map (see e.g. [BPV]), the period domain $\Omega(M^\perp)$ 
parameterizes the marked $M$-polarized K3 surfaces, and  a generic 
point $t \in \Omega(M^\perp)$ comes from a marked $M$-polarized K3 
surface $(X,\tau)$ such that $\tau^{-1}(M)=NS(X)$. 

Mirror symmetry of K3 surfaces is well described for marked 
$M$-polarized K3 surfaces especially when the 
orthogonal lattice $M^\perp$ contains a hyperbolic lattice, i.e.  
$M^\perp=U \oplus \vM$. In this case we have the following  
embedding 
$$
M\oplus U \oplus \vM \subset \Lambda_{\text{K3}},
\Eqno\MvMcond
$$
where $M$ and $\vM$ are of signature $(1,t)$ and $(1,18-t)$ respectively. 
In this paper we always assume the above property ({\MvMcond}) for the lattice 
$M$. And we say that the family of marked $M$-polarized K3 surfaces is mirror 
symmetric to the family of marked $\vM$-polarized K3 surfaces, following [Do].

For the description of mirror symmetry, it will turn out that 
the Mukai lattice is more natural than the K3 lattice 
$\Lambda_{\text{K3}}$. Corresponding to the K3 lattice, let us define 
{\it abstract Mukai lattice} 
$(\tilde \Lambda, \langle \;,\;\rangle)$ by 
$$
\tilde \Lambda = \Lambda^0 \oplus \Lambda_{\text{K3}} \oplus \Lambda^4, 
$$
where $\Lambda^0:=\bold Z e, \Lambda^4:=\bold Z f$ and we naturally 
extend the bilinear form of $\Lambda_{\text{K3}}$ to 
$\tilde \Lambda$ by setting $\langle e,e \rangle=\langle f,f \rangle=0, 
\langle e, f \rangle =-1$, $\langle e, \Lambda_{\text K3} \rangle=\langle f,
\Lambda_{\text{K3}}\rangle=0$. Thus $\tilde H(X,\bold Z)$ is isomorphic 
to $\tilde \Lambda$. 
We call an isomorphism of the lattice 
$\tilde\tau: \tilde H(X,\bold Z) \rightarrow \tilde \Lambda$ a 
{\it Mukai marking}, 
and a Mukai marking $\tilde\tau$ satisfying 
$\tilde\tau(H^0(X,\bold Z))=\Lambda^0, \tilde\tau(H^4(X,\bold Z))=\Lambda^4, 
\tilde\tau(H^2(X,\bold Z))=\Lambda_{\text{K3}}$ a {\it graded Mukai marking}. 
Note that a Mukai marking $\tilde\tau: \tilde H(X,\bold Z) \rightarrow 
\tilde \Lambda$ is not graded in general.  Also note that  
a marking $\tau:H^2(X,\bold Z) \rightarrow \Lambda_{\text{K3}}$ 
naturally extends to a graded Mukai marking $\tilde\tau$ (up to 
$\text{Aut}(\Lambda^0)=\text{Aut}(\Lambda^4)=\{ \pm 1 \}$). 
Now we can state mirror symmetry of the marked K3 surfaces in terms of 
the Mukai lattice. 

\Proclaim{Proposition}{\MirrorMukai} Let $(X,\tau)$ and $(\vX,\vtau)$, 
respectively, be generic marked $M$-polarized K3 surfaces and 
$\vM$-polarized K3 surfaces. Then we have the following identifications 
in the abstract Mukai lattice, 
$$
\overbrace{
\underbrace{
U \;\; \oplus \;\; M}_{\tilde\vtau(T(\vX))} }^{\tilde\tau(\widetilde{NS}(X))}
\;\; \oplus \;\; 
\underbrace{
\overbrace{ U \;\; 
\oplus \;\; \vM }^{\tilde\tau(T(X))}}_{\tilde\vtau(\widetilde{NS}(\vX))}  
\subset \tilde \Lambda
\Eqno\mirrorMukaiL
$$
where $\widetilde{NS}:=H^0\oplus NS \oplus H^4 \cong U\oplus NS$. 
 \qed
\endproclaim

Hereafter we fix the basis $e,f$ for the first hyperbolic lattice 
( the basis we have introduced for the abstract Mukai lattice) and 
the corresponding basis $\ve, \vf$ for the second hyperbolic lattice $U$ 
in (\mirrorMukaiL). 
We assume the intersections to be $\langle \ve,\ve\rangle=
\langle \vf,\vf \rangle=0$ and $\langle \ve,\vf \rangle=-1$.

\noindent
{\bf Remark.} For the lattice $U\oplus \vM$ in the above embedding 
(\mirrorMukaiL), we can associate two different, but isomorphic, domains. 
One is the period domain $\Omega(U\oplus \vM)$ which describes the 
complex structure deformation space of $(X,\tau)$, and the other is a  
tube domain, 
$$
T_{\vM}:=\vM\otimes \bold R + i V(\vM\otimes \bold R)
$$
where $V(\vM\otimes \bold R):=\{ x \in \vM\otimes \bold R \;\vert \; 
(x,x)>0 \}$. 
This tube domain is understood as a covering of the complexified K\"ahler 
moduli space of $(\vX,\vtau)$. 
One can see the mirror correspondence explicitly in the map 
$\mu: T_{\vM} \rightarrow \Omega(U\oplus \vM)$ defined by 
$$
\mu(x)=\bold C ({1\over2}\langle x,x \rangle \vf + x + \ve)  \;,\;\; 
(x=\vB+i \vK \in T_{\vM}), 
\Eqno\TtoOmega
$$  
which is called the {\it mirror map}. We can verify easily that this map 
is bijective using a property 
$\langle \omega, \vf \rangle \not=0$ for $\bold C \omega 
\in \Omega(U\oplus\vM)$( see [Do], Lemma (4.1)).   \qed

\vskip0.5cm

Let $\tilde \tau: \tilde H(X,\bold Z) \rightarrow \tilde \Lambda$ be a 
Mukai marking. We consider the group 
$$ 
O_{Hodge}(\tilde \Lambda,\tilde\tau(\bold C \omega_X))
:=\{ g \in O(\tilde \Lambda) \vert 
g(\tilde\tau(\bold C \omega_X))=\tilde \tau(\bold C \omega_X) \}.
$$
Then $O_{Hodge}(\tilde \Lambda,\tau(\bold C \omega_X))\cong O_{Hodge}(\tilde 
H(X,\bold Z))$ (cf. Definition {\HodgeIso}). 

Now we can argue two different pictures for the same group 
$O_{Hodge}(\tilde \Lambda,\tilde\tau(\bold C \omega_X))$ based on the 
mirror relation (\mirrorMukaiL). 

The first is to understand this group as the Hodge isometries of 
marked $M$-polarized K3 surface $(X,\tau)$ and its 
extension $(X,\tilde\tau)$ to a graded Mukai marking. By our Theorem 
{\MainThI}, this group contains the image $ch(\text{Auteq} D(X))$ 
as a subgroup of index at most two. 

The second is the mirror picture to the first and valid for the mirror 
$(\vX,\vtau)$ and its extension $(\vX,\tilde\vtau)$ to 
a graded Mukai marking. To describe this 
let us recall the mirror relation $\tilde\tau(T(X)) = U\oplus \vM = 
\tilde\vtau(\widetilde{NS}(\vX))$ for  generic $(X,\tilde\tau)$ and 
$(\vX,\tilde\vtau)$. 
Among the Hodge isometries in  
$O_{Hodge}(\tilde\Lambda,\tilde\tau(\bold C \omega_X))$, 
let us focus on the isometries of $\tilde\Lambda$ which stabilize  
$\tilde\tau(T(X))$  and act as identity on $\tilde\tau(T(X))$. 
We denote the subgroup consisting of these isometries by 
$$
O^1_{Hodge}(X,\tilde\tau):=\{ 
\varphi \in O_{Hodge}(\tilde\Lambda,\tilde\tau(\bold C \omega_X)) 
\; \vert \; \varphi|_{U\oplus \vM}=
\text{id}_{U\oplus \vM} \}.
\Eqno{\Oid}
$$
Now consider an index two subgroup $O^{1,+}_{Hodge}
(X,\tilde\tau)$ of $O^1_{Hodge}(X,\tilde\tau)$ 
which preserves the orientations of the 
positive two planes in $(U\oplus M)\otimes\bold R$. 
Recall the isomorphism $\mu:T_{\vM} \simrightarrow \Omega(U\oplus \vM)$ 
(\TtoOmega) and observe that elements in 
$O^{1,+}_{Hodge}(X,\tilde\tau)$ 
preserve each complexified K\"ahler class of the mirror $(\vX,\tilde\vtau)$  
defined by $\vB+i\vK:=\mu^{-1}(\bold C \omega_X)$. 
(On the mirror side, presumably $ch(\text{Auteq} D(X))=O_{Hodge}^+(\tilde 
H(X,\bold Z))$ should hold, where $O_{Hodge}^+(\tilde H(X,\bold Z))$ is the 
index two subgroup which preserves the orientations of the positive 
four planes in $\tilde H(X,\bold Z)\otimes \bold R$. 
(See [Sz, Conjecture 5.4].) 
Note that the tube domain 
$T_{\vM}$ and the period domain $\Omega(U\oplus \vM)$ both have 
two connected components. This also seems closely related to such 
index-two phenomena in the present paper.)

Now recall that the kernel  
$O(T(\vX))^*\cong O(\tilde\vtau(T(\vX)))^*$ is a subgroup of $O(T(\vX))$ 
consisting of those elements $\phi_{T(\vX)}$ such that 
$(\phi_{T(\vX)}, \text{id}_{\widetilde{NS}(\vX)})$ $\in 
O(T(\vX))\times O(\widetilde{NS}(\vX))$ extend to elements of 
$O(\tilde H(\vX,\bold Z))\cong O(\tilde \Lambda)$ (cf. Remark after 
Theorem {\MainThII}). Using this extension property, we can identify 
the subgroup $O^{1}_{Hodge}(X,\tilde\tau)$ with $O(T(\vX))^*$ and 
further its index two subgroup $O^{1,+}_{Hodge}(X,\tilde\tau)$ with 
$O^+(T(\vX))^*$. 

Now we can apply Theorem {\MainThII} for $Y=\vX$ to see that 
all elements in $O^{1,+}_{Hodge}(X,\tilde\tau)$ 
come from the symplectic mapping class group of $(\vX,\vtau)$ 
with respect to its K\"ahler class $\text{Im}(\mu^{-1}(\bold C \omega_X))$.  
This is the mirror interpretation of the subgroup 
$O_{Hodge}^{1,+}(X,\tilde\tau)$ $\subset$  
$O_{Hodge}(\tilde\Lambda,\tilde\tau(\bold C \omega_X))$.  
Here we have arrived at a subgroup $O_{Hodge}^{1,+}(X,\tilde\tau)$ 
which has the index greater than two in 
$O_{Hodge}(\tilde\Lambda,\tilde\tau(\bold C \omega_X))$. 
This should be attributed to the fact that we have replaced the desired 
autoequivalence group of ``$D Fuk(\vX)$'' by the well-known but 
possibly smaller  symplectic mapping class group. 
Note, for example, that the shift functor 
contained in $\text{Auteq} D(X)$ does not have its counter part 
in the symplectic mapping class group  but should have 
in ``$\text{Auteq} D Fuk(\vX)$''.

\vskip0.5cm
\noindent
(1-4) Fourier-Mukai partners and monodromy of the mirror 
family 

\vskip0.2cm

When a derived category $D(X)$ is given, we can ask for the varieties $Y$ 
which admit the equivalence $\Phi: D(X) \cong D(Y)$. The smooth projective 
varieties with this 
property are called {\it Fourier-Mukai(FM) partners} of $X$. If $X$ has ample 
canonical or anticanonical bundle, Bondal and Orlov[BO] proved that $X$ itself 
is the only FM partner.  For K3 surfaces, however, $X$ has in general finitely 
many FM partners:

\Proclaim{Proposition}{\BMfiniteness} {\bf ([BM, Proposition 5.3]}, see 
also {\bf [Og2, Proposition ({\preMainII})])} For a given K3 surface $X$, 
there are only finitely many FM partners.   \qed
\endproclaim

For a generic K3 surface, or more precisely, for a 
K3 surfaces of $\rho(X)=1$, the number of FM partners has been determined 
as follows;

\eject

\Proclaim{Proposition}{\FMnumber} {\bf ([Og2, Proposition (1.10)])} 
Let $X$ be a K3 surface with $NS(X)=\bold Z h$. Set $\text{deg}(X)=(h^2)
=2 n$. Then the number of FM partners of $X$ is given by $2^{p(n)-1}$, where 
$p(1)=1$ and $p(n) \; (n\geq2)$ is the number of prime numbers $p \; (\geq 2)$ 
such that $p | n$.  \qed
\endproclaim

We will arrive at the same number studying the monodromy representation 
of the mirror family $\vX$. To define the mirror family let us first 
remark that a lattice of rank one admits a unique 
primitive embedding into the K3 lattice $\Lambda_{\text{K3}}$. Now let us 
consider a rank one lattice $M_{n}=\langle 2n \rangle$, i.e.  
a lattice $M_n=\bold Z v$ with its bilinear form determined by $(v,v)=2n$. 
Then 
because of the uniqueness (up to isomorphism) of the primitive embedding 
we have the following decomposition;
$$
M_n \oplus U \oplus \vM_n  
\subset \Lambda_{\text{K3}},  
$$
where $\vM_n:=\langle -2n \rangle \oplus U 
\oplus E_8(-1)^{\oplus2}$. 
A K3 surface with $NS(X)=\bold Z h$ and $\text{deg}(X)=2n$ may be regarded 
as a generic member of the family of the marked $M_n$-polarized K3 surfaces. 
The mirror family defined in (1-3) is the marked $\vM_n$-polarized K3 surfaces 
$(\vX,\vtau)$, which are parametrized by the period domain 
$\Omega(T(\vX))$ $=\Omega(U\oplus M_n)$. The generic member 
$(\vX,\vtau)$ of the 
family has its transcendental lattice $T(\vX) \cong U \oplus M_n$ and 
{\it classified} by $\Omega^0(U\oplus M_n)$, a complement of countable union 
of proper closed subsets of $\Omega(U\oplus M_n)$. 
($\Omega(U\oplus M_n)\setminus \Omega^0(U\oplus M_n)$ is however 
dense in $\Omega(U\oplus M_n)$, see e.g. [Og2]). 
For generic elements of the marked 
$\vM$-polarized K3 surfaces, we have: 

\Proclaim{Lemma}{\genericLemma} Let $(\vX_1,\vtau_1)$ and $(\vX_2, \vtau_2)$ 
be marked $\vM_n$-polarized K3 surfaces pa-rametrized by 
$\Omega^0(U\oplus M_n)$. Then $\vX_1 \cong \vX_2$ if and only if 
there exists $g \in O(U\oplus M_n)$ such that $g(\bold C \tau_1(\omega_{X_1}))=
\bold C \tau_2(\omega_{X_2})$. 
\endproclaim

The proof of this lemma will be given in the section 4. Using the natural 
action $\iota: O(U\oplus M_n) \rightarrow \text{Aut} \Omega^0(U\oplus M_n)$, 
we consider the quotient space $\Omega^0(U\oplus M_n)/O(U\oplus M_n) 
\cong \Omega^{0,+}(U\oplus M_n)/O^+(U\oplus M_n)$, where 
$\Omega^{0,+}(U\oplus M_n)=\Omega^0(U\oplus M_n)\cap \Omega^+(U\oplus M_n)$ 
and $\Omega^{+}(U\oplus M_n)$ is one of the two connected components of 
$\Omega(U\oplus M_n)$. 
This space is the classifying space of the generic $\vM_n$-polarizable 
K3 surface by Lemma {\genericLemma}. 
The closure 
of $\Omega^{0,+}(U\oplus M_n)/O^+(U\oplus M_n)$ is 
$\Omega^+(U\oplus M_n)/O^+(U\oplus M_n)$. 
Note that $\pm \text{id}_{U\oplus M_n}$ 
acts trivially on the period domain, i.e.  $\text{Ker}(\iota)=
\{\pm \text{id}_{U\oplus M_n} \}$. This leads to the following definition:

\Proclaim{Definition}{\monodG} We call the group $O^+(U\oplus M_n)/\{ \pm 
\text{id}_{U\oplus M_n} \}$ the monodromy group of the 
$\vM_n$-polarizable K3 surfaces $\vX$, and denote it by ${\Cal M}_n(\vX)$. 
\qed
\endproclaim

Now consider the definition (\Ostar) for $T(Y)=T(\vX)\cong U\oplus M_n$, and 
see the following composition of natural maps;
$$
O^+(U\oplus M_n)^* \rightarrow O^+(U\oplus M_n) 
\rightarrow O^+(U\oplus M_n)/\{ 
\pm \text{id}_{U\oplus M_n} \} =:{\Cal M}_n(\vX).
\Eqno{\compOp}
$$
Note that $-\text{id}_{U\oplus M_n}$ is contained in $O^+(U\oplus M_n)^*$ 
only for $n=1$ since $A((U\oplus M_n)^*/(U\oplus M_n))=A(M_n^*/M_n)=
\langle v/2n \rangle \cong \bold Z/2n$. Therefore we see that 

\Proclaim{Lemma}{\compo} The composition map $O^+(U\oplus M_n)^* 
\rightarrow {\Cal M}_n(\vX)$ defined in (\compOp) is injective for $n\geq 2$ 
and has the kernel $\{ \pm \text{id}_{U\oplus M_n} \}$ for $n=1$. 
\qed
\endproclaim

Applying our Theorem {\MainThII} to $Y=\vX$ and using $T(Y)\cong 
U\oplus M_n$ for a generic $\vM_n$-polarizable K3 surface, we have
$$
O^+(U\oplus M_n)^* = \text{Im}(\pi_0\text{Symp}(\vX,\kappa_{\vX}) \rightarrow 
O(U\oplus M_n)). 
$$
Based on this relation and Lemma {\genericLemma}, we define:

\Proclaim{Definition}{\sympMonod} For a generic symplectic structure 
$(\vX, \kappa_{\vX})$ of a generic $\vM_n$-polarizable K3 surface $\vX$, 
we call the group $O^+(U\oplus M_n)^*$ for $n\geq2$ ( respectively, the group  
$O^+(U\oplus M_n)^*/\{ \pm \text{id}_{U\oplus M_n} \}$ for 
$n=1$) the monodromy representation of the symplectic mapping class  
group of $(\vX, \kappa_{\vX})$. We denote this group by $\Cal{MS}_n(\vX)$.
\qed
\endproclaim

Now we can state our theorem:

\Proclaim{Theorem}{\MainThIII} {\bf (Main Theorem 3)} Let $X$ be a K3 surface 
of $\rho(X)=1$ and $\text{deg}(X)=2n$. Then the number of FM partners 
of $X$ is given by the index $[{\Cal M}_n(\vX):\Cal{MS}_n(\vX)]$. 
\endproclaim

Recall that the number of FM partners is $2^{p(n)-1}$ by  
Proposition \FMnumber. Our theorem above derives the same numbers 
from the monodromy property of the mirror family $\vX$. A proof of Theorem 
{\MainThIII} will be given in section 4. 

In section 5, we will also study in details the first non-trivial 
case of $n=6$. We will present the monodromy calculations explicitly 
following [LY2][PS][BP] 
and show how the monodromy property is connected to the numbers of 
FM partners in this particular case.

\head
{\S 2 Autoequivalences and Proof of Theorem \MainThI}
\endhead

\global\secno=2  
\global\propno=1 
\global\eqnum=1 

\noindent
(2-1) Various autoequivalences.  
\vskip0.2cm

For our proof of Theorem {\MainThI}  
let us recall basic autoequivalences in order. 

\noindent
1) {\it Shift functor $[n]$:} $D(X) \rightarrow D(X) (n \in \bold Z)$ 
defined by $K^\bullet \rightarrow L^\bullet=K^{\bullet+n}$, i.e.  the shift 
by $n$ to the left, is an autoequivalence. This functor does not change 
the complex except its order-preserved numbering. However we should 
note that $f^{[n]}(:=ch([n]))=-\text{id}$ if $n\equiv 1 \; (2)$ and 
$f^{[n]}=\text{id}$ if $n\equiv 0 \; (2)$ by the definition $ch(K^\bullet)
=\sum_i (-1)^i ch(K^i)$. 

\noindent
2) $Aut(X)$: An automorphism $g \in \text{Aut}(X)$ 
gives rise to an autoequivalence  
$g: D(X) \rightarrow D(X)$ by sending $K^\bullet$ to $L^\bullet$ 
with $L^i = g^* K^i$. Since the quasi-inverse to $g$ is given by $g^{-1}$, 
$g$ is an autoequivalence. $g$ is represented by a kernel 
$\Cal O_{\Gamma(g)} \in D(X \times X)$, i.e.  $g=\Phi^{\Cal O_{\Gamma(g)}}$ 
as an element of $\text{Auteq} D(X)$, where 
$\Gamma(g)=\{ (x,g(x)) | x \in X \} 
\subset X \times X$ is the graph of $g$ and $\Cal O_{\Gamma(g)}$ is the 
structure sheaf of the reduced closed subscheme $\Gamma(g) \subset X \times X$ 
( which we identify the pushforward $\iota_* \Cal O_{\Gamma(g)}$ under  
$\iota :\Gamma(g) \hookrightarrow X \times X$). The induced action 
$f^{\Cal O_{\Gamma(g)}}$ on $\tilde H(X,\bold Z)$ is the pullback 
by $g^*$. 

\noindent
3) {\it Tensoring by line bundles:} Let $\Cal L \in \text{Pic}(X)$ be a line 
bundle (invertible sheaf) on $X$. Then we may associate to it an 
autoequivalence $\Phi^{\Cal L}:=\Phi^{\pi_2^* \Cal L}: 
D(X) \rightarrow D(X)$ by $\Cal X \mapsto \pi_{2*}
(\pi_2^* \Cal L \otimes \pi_1^* \Cal X)=\Cal L \otimes \pi_{2*}
(\pi_1^* \Cal X)$, where $\pi_1$ and $\pi_2$ are the natural projections 
$\pi_{1,2}:X \times X \rightarrow X$ to the first and the second $X$, 
respectively. 
The quasi-inverse of this is simply given by $\Phi^{\Cal L^{-1}}$. The induced 
action $f^{\Cal L}:=f^{\pi_2^*\Cal L}$ on $\tilde H(X,\bold Z)$ is the  
multiplication by the chern character $ch(\Cal L)=(1, c_1(\Cal L), {1\over 2} 
c_1(\Cal L)^2)$ in the graded ring $\tilde H(X,\bold Z)$. 

These three functors 1),2),3) above form a subgroup of $\text{Auteq} D(X)$ 
which is isomorphic to $(\bold Z \times \text{Pic} X) \rtimes Aut(X)$. 

\noindent
4) {\it Twistings by spherical objects [ST]:} Let $C$ be a smooth rational 
curve in a K3 surface $X$, then $(C^2)=-2$.  In this paper, we mean 
by smooth $(-2)$ curve a smooth rational curve. Now consider 
the structure sheaf $\Cal O_C(:=\iota_* \Cal O_C)$. Using the exact sequence 
$0 \rightarrow \Cal O_X(-C) \rightarrow \Cal O_X \rightarrow \Cal O_C 
\rightarrow 0$ and evaluating $\text{Hom}(*, \Cal O_C)$, we have 
$\text{Ext}^i(\Cal O_C, \Cal O_C)= \bold C$ for $i=0,2$ and 
$\text{Ext}^i(\Cal O_C, \Cal O_C)= 0 \; (i\not=0,2)$. This a simple 
example what is called a {\it spherical object} in the derived category[ST]. 
Here we simply summarize relevant results restricting our attentions to 
the case of K3 surface $X$. 

\Proclaim{Theorem}{\ST} {\bf ([ST, Theorem 1.2])} Let $X$ be a K3 surface. 
$\Cal E \in D(X)$ is called spherical if it satisfies
$$
\text{Hom}^i_{D(X)}(\Cal E,\Cal E)=
\cases \bold C & i=0,2 \cr
       0       & i\not=0,2 \;\;.\cr
\endcases
$$
(Note that $\Cal O_X(K_X)\cong \Cal O_X$.) For a spherical object $\Cal E$ we 
consider the mapping cone 
$\Cal C:=\text{Cone}(\pi_1^* \Cal E^\vee \otimes \pi_2^* \Cal E 
\rightarrow \Cal O_\Delta)$ of the natural evaluation 
$ \pi_1^* \Cal E^\vee \otimes \pi_2^* \Cal E \rightarrow \Cal O_\Delta$,  
where $\Delta \subset X\times X$ is the diagonal 
and $\Cal O_\Delta$ is at the zeroth position of the complex $\Cal C$. 
$\Cal E^\vee = \bold R \Cal Hom(\Cal E,\Cal O_{X\times X})$ is the derived 
dual of $\Cal E$, and $\pi_1$ and $\pi_2$ are the natural projections 
$\pi_{1,2}: X\times X 
\rightarrow X$ to the first and the second, respectively. Then the functor 
$T_{\Cal E}:=\Phi^{\Cal C}:D(X)\rightarrow D(X)$ 
defines an equivalence. This functor is called a twist functor. 
The corresponding map $t_{\Cal E}:=f^{\Cal C}$ on $\tilde H(X,\bold Z)$ 
is given by 
$$
t_{\Cal E}(x)=x+\langle ch(\Cal E)\sqrt{td_X}, x \rangle ch(\Cal E)\sqrt{td_X}.
\Eqno{\twistCh}
$$
\endproclaim

\Proclaim{Lemma}{\reflection} Let $T_{\Cal E} \;(\Cal E=\Cal O_C)$ 
be the twist functor with respect to a smooth $(-2)$ curve $C$ 
in a K3 surface $X$ and consider an autoequivalence 
$\Phi^{\Cal O_X(C)}\circ T_{\Cal E}: D(X) 
\rightarrow D(X)$. Then the corresponding action $f^{\Cal O_X(C)}\circ 
t_{\Cal E}$ in $\tilde H(X,\bold Z)$ is a Hodge isometry of 
$(\tilde H(X,\bold Z), \langle\;,\;\rangle)$ and coincides with 
the reflection by the curve $C$,
$$
f^{\Cal O_X(C)}\circ t_{\Cal E}(x)=x+\langle x, C\rangle C .
$$
\endproclaim

\demo{Proof} From the exact sequence $0\rightarrow \Cal O_X(-C) \rightarrow 
\Cal O_X \rightarrow \Cal O_C \rightarrow 0$, we have $ch(\Cal O_C)=
ch(\Cal O_X)-ch(\Cal O_X(-C))=(1,0,0)-(1,-C,{1\over2}(-C)^2)=
(0,C,1)$. Now for $x=(\alpha,\beta,\gamma) \in \tilde H(X,\bold Z)=
H^0(X,\bold Z)\oplus H^2(X,\bold Z)\oplus H^4(X,\bold Z)$ we 
have
$$
\aligned
t_{\Cal E}(x)
&=x+\langle ch(\Cal O_C) \sqrt{td_X}, x \rangle ch(\Cal O_C) \sqrt{td_X}\cr
&= x+\langle (0,C,1)\cdot (1,0,1), (\alpha,\beta,\gamma)\cdot (1,0,1) \rangle
  (0,0,1) \cr
&=(\alpha,\beta,\gamma)+(-\alpha+(\beta,C))\times(0,C,1) \cr 
&=(\alpha,\beta+(-\alpha+(\beta,C))C,\gamma-\alpha+(\beta,C)) \cr
\endaligned
$$
Now we note that $f^{\Cal O_X(C)}(x)=ch(\Cal O_X(C))x$, i.e.  simply 
the multiplication by $ch(\Cal O_X(C)) =(1,C,{1\over2}C^2)$. 
Composing these two actions we obtain the desired result,
$$
f^{\Cal O_X(C)}\circ t_{\Cal E}(x)
=(\alpha,\beta+(\beta,C) C, \gamma) 
= x+\langle x,C \rangle C 
\qed
$$
\enddemo

\noindent
5) {\it Switching functor (a special case of the Fourier-Mukai transform 
for fine moduli spaces of stable sheaves):} Let us first recall the following 
result due to Mukai. 

\Proclaim{Theorem}{\MukaiModuli} 
{\bf ([Mu2]} see also {\bf [BM, Cor.2.8])} 
Let $X$ be a K3 surface with a fixed polarization. Consider a smooth 
fine compact two dimensional moduli space $Y$ of stable sheaves on $X$ 
and denote by $\Cal P$ the universal sheaf on $X\times Y$ ( for which 
we have for each $ y \in Y$ the stable sheaf $\Cal P_y=\pi_Y^*\Cal O_y \otimes 
\Cal P$ on $X\times\{y\}\cong X$, which represents the point $y$ in the 
moduli space $Y$).  
Then $Y$ is smooth, hence a K3 surface and 
$\Phi_{Y\rightarrow X}^{\Cal P}:D(Y)\rightarrow D(X)$ 
is a FM transform and satisfies $f^{\Cal P}(ch(\Cal O_y))=ch(\pi_{X*}
\Cal P_y)\sqrt{td_X}$. (Note that $ch(\Cal O_y)\sqrt{td_X}=ch(\Cal O_y)=
(0,0,1)$.)  \qed 
\endproclaim

\noindent
For our purpose we apply this theorem to a special case. Let $\Delta \subset 
X \times X$ be the diagonal and $\Cal I_\Delta$ be its ideal sheaf in 
$X \times X$. We can regard $X$ as the fine moduli space of the 
ideal sheaves of points $\Cal I_x (x \in X)$, which are certainly stable. 
Therefore $\Cal I_\Delta$ is the universal sheaf on $X \times X$. 
By the above theorem we have the 
corresponding FM transform (autoequivalence) $\Phi^{\Cal I_\Delta}:
D(X) \rightarrow D(X)$ for which we have $f^{\Cal I_\Delta}(\alpha,\beta,
\gamma)=(\gamma,-\beta,\alpha)$. The last equation for $f^{\Cal I_\Delta}$ 
follows from $ch(\Cal I_\Delta)=ch(\Cal O_{X\times X})-
ch(\Cal O_\Delta)$ and $f^{\Cal I_\Delta}(x)=-x-\langle  
ch(\Cal O_X) \sqrt{td_X}, x \rangle ch(\Cal O_X) \sqrt{td_X}$, i.e. 
$$
f^{\Cal I_\Delta}(\alpha,\beta,\gamma)=-(\alpha,\beta,\gamma) 
- \langle (\alpha,\beta,\gamma),(1,0,1) \rangle (1,0,1) 
=(\gamma,-\beta,\alpha).
$$
We call this autoequivalence $\Phi^{\Cal I_\Delta}$ a 
{\it switching functor}.

\vskip0.5cm
\noindent
(2-2) FM transforms on a K3 surface. 
\vskip0.2cm

Here we recall a theorem of Mukai and Orlov on  FM transforms on a K3 surface: 

\Proclaim{Theorem}{\FM-Hodge} {\bf ([Mu1,2], [Or1],} see also {\bf [BM] )} 
Let $X$ and $Y$ be K3 surfaces. Then the 
following statements are equivalent;
\item{1)} there exists a FM transform $\Phi: D(Y) \rightarrow D(X)$
\item{2)} there exists a Hodge isometry $f_{T}: (T(Y),\bold C \omega_Y) 
\rightarrow (T(X),\bold \omega_X)$
\item{3)} there exists a Hodge isometry $f: (\tilde H(Y,\bold Z), \bold C 
\omega_Y) \rightarrow (\tilde H(X,\bold Z), \bold C \omega_X)$ 
\item{4)} $Y$ is isomorphic to a two dimensional fine moduli space 
of stable sheaves on $X$.
\endproclaim

In the following arguments, since we need not only the results of 
Theorem {\FM-Hodge} but also the arguments, we sketch proof here along [BM]: 

\noindent 
1)$\Rightarrow$2) Writing the FM transform 
$\Phi=\Phi^{\Cal E}_{Y\rightarrow X}$ by a 
kernel $\Cal E$, one has a Hodge isometry $f^{\Cal E}_{Y\rightarrow X}: 
(\tilde H(Y,\bold Z),\bold C \omega_Y)\rightarrow 
(\tilde H(X,\bold Z),\bold C \omega_X)$ 
(Theorem {\D(X)}). Recall that the transcendental lattice is the minimal 
primitive sublattice of $\tilde H(X,\bold Z)$ which satisfies 
$\bold C \omega_X \subset T(X) \otimes \bold C$.  
Therefore the restriction $f^{\Cal E}_{Y\rightarrow X}$ to $T(Y)$ gives 
the desired 
Hodge isometry $f_T:=f^{\Cal E}_{Y\rightarrow X}\vert_{T_Y}:
(T(Y),\bold C \omega_Y) \rightarrow (T(X),\bold C \omega_X)$. 

\noindent
2)$\Rightarrow$3) By Nikulin's theorem of primitive embedding of lattices 
[Ni, Theorem (1.14.4)] $f_T:
(T(Y), \bold C \omega_Y)\rightarrow (T(X), \bold C \omega_X)$ extends to $f:
(\tilde H(Y,\bold Z), \bold C \omega_Y) \rightarrow 
(\tilde H(X,\bold Z), \bold C \omega_X)$ satisfying $f\vert_{T(Y)}=f_T$. 

 \noindent
3)$\Rightarrow$4) Since this part is crucial for our purpose, we present 
this in detail. Let $f: (\tilde H(Y,\bold Z),\bold C \omega_Y) 
\rightarrow (\tilde H(X,\bold Z),\bold C \omega_X)$ be a Hodge isometry and 
set $v=f((0,0,1))$. Let us first show that, {\it 
if necessary composing $f$ with 
a suitable $f^{\Cal E}_X=ch(\Phi_X^{\Cal E})$ of 
$\Phi^{\Cal E}_X \in \text{Auteq} D(X)$, one may assume 
for $v=(r,l,s)$ that $r>1$, $l$ is ample and $s$ is coprime to $r$. } 
First of all let us note that, since $l$ in $(r,l,s)=f((0,0,1))$ is algebraic, 
we may assume $r>1$ considering a composition 
$f^{\Cal L}\circ f$ with a suitable line bundle and further with  
the Hodge isometries $f^{\Cal I_\Delta}, f^{[1]}=- \text{id}$ 
associated to the switching functor $\Phi^{\Cal I_\Delta}$ and 
the shift functor $[1]$, respectively. We redefine $(r,l,s)$ to 
be the image of the vector $(0,0,1)$ under the compositions to ensure 
$r>1$.  Now we show that one can satisfy the condition $(r,s)=1$ further  
compositions of suitable Hodge isometries. For this purpose let us write 
$u=(a,b,c)=f((1,0,0))$. Since $f$ is a Hodge isometry we have 
$\langle u,v \rangle = \langle (1,0,0),(0,0,1)\rangle=-1$. 
On the other hand, by the definition of $\langle \;,\; \rangle$, we 
calculate $\langle u,v \rangle=-as-cr+(b,l)$. Therefore we have 
$-1= -as-cr+(b,l)$, which means $s,r,(b,l)$ 
are coprime. The last condition ensures that there exists an integer 
$n$ for which $r$ and $s+n (b,l)$ are coprime. 
Now consider a composition $
f^{\Cal O_X(n b)}\circ f$ for which 
we have 
$$
f^{\Cal O_X(n b)}\circ f((0,0,1))=f^{\Cal O_X(n b)}(v)=
(r,l+r n b,s+n(b,l)+{n^2 r\over2}(b,b)). 
$$ 
Since $r$ and $s+n(b,l)$ are 
coprime, redefining $(r,l,s)$ to be $f^{\Cal O_X(b)}\circ f((0,0,1))$ 
we have $(r,s)=1$ and $r>1$. Finally consider a composition by 
$f^{\Cal O_X(r A)}$ 
with $A$ being sufficiently ample. Then we have 
$$
f^{\Cal O_X(r A)}(r,l,s)=(1,rA,{1\over2}(A,A))(r,l,s)=
(r,r^2 A+l,{1\over2}r^3 (A,A)+s),
$$ 
i.e. we obtain an ample class for the second factor preserving the other 
conditions. 

Since all the functors we used so far  
are Hodge isometries, we may assume that $v=f((0,0,1))=(r,l,s)$ has 
the desired properties from the beginning. 

Now since $(0,0,1)$ is algebraic and $f$ is a Hodge isometry, $v=f((0,0,1))$ 
is also algebraic and lies in $H^0(X,\bold Z)\oplus 
NS(X) \oplus H^4(X,\bold Z)$, which is perpendicular 
to $T(X)$. Therefore one can consider the moduli space $Y^+$ of stable 
sheaves on $X$ whose Mukai vector is $v$, i.e.  stable sheaves $\Cal E$ with 
$ch(\Cal E)\sqrt{td_X}=v$ with respect to the ample polarization $l$. Recall 
that since $v$ is primitive $Y^+$ is fine [Mu2, Theorem A.6] and $\langle 
v, v \rangle=\langle (1,0,0),(1,0,0)\rangle=0$ implies $\text{dim} Y^+=2$. 
Moreover $Y^+$ is smooth, non-empty and compact. 
This smoothness follows form the Main Theorem 
in [Mu1], non-emptiness follows from Theorem 5.4 in [Mu2] and 
compactness follows from Proposition 4.1 in [Mu2] and $(r,s)=1$. Then by the 
main Theorem of Bridgeland [Br, Theorems 5.1 and 5.3] (see also [Or1] for 
another argument) one has an equivalence $\Phi^{\Cal E}_{Y^+ \rightarrow 
X}: D(Y^+) \rightarrow 
D(X)$ for which $f^{\Cal E}_{Y^+\rightarrow X}((0,0,1))=v$. Then, 
using 1)$\Rightarrow$3), $g:=(f^{\Cal E})^{-1}
\circ f: \tilde H(Y,\bold Z) \rightarrow \tilde H(Y^+,\bold Z)$ is a 
Hodge isometry satisfying $g((0,0,1))=(0,0,1)$. This Hodge isometry $g$ 
of the Mukai lattices reduces to that of the second cohomologies $H^2$, i.e. 
$g: H^2(Y,\bold Z)=
(0,0,1)^\perp/\bold Z (0,0,1) \simrightarrow H^2(Y^+,\bold Z)=
(0,0,1)^\perp/\bold Z (0,0,1)$. Therefore by the Torelli Theorem, we conclude 
$Y \cong Y^+$. 

\noindent
4)$\Rightarrow$1) This is a special case of the Theorem {\MukaiModuli}.  
\qed

\vskip0.5cm
\noindent
(2-3) Proofs of Theorem {\MainThI} 
\vskip0.2cm

\noindent
{\it Proof of Theorem {\MainThI}:} 
Now we come to a proof of our main Theorem. Let us consider a K3 surface 
$X$ with a fixed {\it graded} Mukai markings 
$\tilde\tau_X: \tilde H(X,\bold Z) \rightarrow 
\tilde \Lambda$. We set $\omega_X':= \tilde\tau_X(\omega_X)$. Then 
let $g\in O_{Hodge}(\tilde H(X,\bold Z))=O_{Hodge}(\tilde \Lambda, \bold C 
\omega_X')$ and consider $g((0,0,1))$. 
As in the proof of Theorem \FM-Hodge we consider   
compositions of $g$ with the Hodge isometries 
$f^{\Cal L}_{X}, f^{\Cal I_\Delta}_{X}$, 
which correspond to autoequivalences 
$\Phi_{X}^{\Cal L}$ and $\Phi^{\Cal I_\Delta}_{X}$, 
respectively. Choosing a suitable composition, we may assume  
$$
F(f^{\Cal L_1}_{X }, f^{\Cal L_2}_{X }, 
f^{\Cal I_\Delta}_{X})\circ g ((0,0,1))=(a,b,c)
$$ 
with 
$a$ and $c$ coprime and $b$ ample, i.e.  $b=\tilde\tau_X(B)$ for an ample 
line bundle on $X$, where $F(f^{\Cal L_1}_{X }, \cdots )$ 
represents the composition we chose. Then, since the marking 
$\tilde\tau_X$ is  graded, and by the argument for 3)$\Rightarrow$ 4) 
of Theorem \FM-Hodge, 
there is a K3 surface $Y$ and $\Cal E \in D(Y\times X)$ such that 
$\Phi^{\Cal E}_{Y\rightarrow X }: D(Y)\rightarrow D(X)$ is an equivalence 
and $f^{\Cal E}_{Y \rightarrow X }((0,0,1))=(a,b,c)$. 
Here the Hodge isometry 
$f^{\Cal E}_{Y\rightarrow X}:(\tilde H(Y,\bold Z),\bold C \omega_Y) 
\rightarrow (\tilde H(X,\bold Z),\bold C \omega_X)$ is identified with the 
Hodge isometry $\tilde\tau_{X}\circ f_{Y\rightarrow X}^{\Cal E}\circ 
\tilde\tau_Y^{-1}: (\tilde \Lambda,\bold C \omega_Y') 
\rightarrow (\tilde \Lambda, 
\bold C \omega_X')$ under the {\it graded } Mukai marking 
$\tilde\tau_Y:\tilde H(Y,\bold Z) \rightarrow \tilde \Lambda$ with $\omega_Y'=
\tilde\tau_Y(\omega_Y)$. 
Now consider a Hodge isometry under this identification 
$$
h:=
(F(f^{\Cal L_1}_{X}, f^{\Cal L_2}_{X}, 
f^{\Cal I_\Delta}_{X})\circ g)^{-1}\circ 
f_{Y\rightarrow X}^{\Cal E}: (\tilde \Lambda, \bold C \omega_Y') 
\rightarrow (\tilde \Lambda, \bold C \omega_X') .
\Eqno\Hodgeh
$$
Then this has the property $h((0,0,1))=(0,0,1)$. Therefore $h$ induces a 
Hodge isometry in the K3 lattice $\Lambda^2=(0,0,1)^\perp/\bold Z (0,0,1)$;
$$
\bar h: ( (0,0,1)^\perp/\bold Z (0,0,1), \bold C \omega_Y') 
\simrightarrow 
( (0,0,1)^\perp/\bold Z (0,0,1), \bold C \omega_X') .
$$
Since the Mukai marking $\tilde\tau_X$ and $\tilde\tau_Y$ are graded, 
$\tilde\tau_X\circ \bar h \circ \tilde\tau^{-1}\vert_{H^2(Y,\bold Z)}$ 
is a Hodge isometry from $(H^2(Y,\bold Z),\bold C \omega_Y)$ to $
(H^2(X,\bold Z),\bold C \omega_X)$. 
This implies that $Y\cong X$ and that 
the map $h$ in (\Hodgeh) is a Hodge isometry from 
$(\tilde \Lambda,\bold C \omega_X')$ to $(\tilde \Lambda,\bold C \omega_X')$. 
Now set $h((1,0,0))=(a,b,c)$. Then since $h$ is an isometry, one has
$$
\aligned
b^2-2ac &=\langle (a,b,c),(a,b,c)\rangle=\langle (1,0,0),(1,0,0)\rangle =0 \cr
-a & =\langle (a,b,c),(0,0,1)\rangle=\langle (1,0,0),(0,0,1)\rangle =-1. \cr
\endaligned
$$
Therefore one knows $h((1,0,0))=(a,b,c)=(1,b,{1\over2}b^2)$. Since 
$(1,0,0)$ is algebraic and $h$ is a Hodge isometry, $(1,b,{1\over2}b^2)$ is 
algebraic, i.e.  $b$ comes from a line bundle, say, $B$. Let  
$f^{(-B)}_X$ be the Hodge isometry corresponding to the autoequivalence 
$\Phi^{(-B)}_{X}$. Then  we have a Hodge isometry
$$
k:=f^{(-B)}_X\circ h: (\tilde \Lambda, \bold C \omega_X') \rightarrow (\tilde 
\Lambda, \bold C \omega_X') 
$$
which satisfies $k((1,0,0))=(1,0,0)$ and $k((0,0,1))=(0,0,1)\cdot(1,-b,
{1\over2}b^2)=(0,0,1)$. In particular we have $k\vert_{\Lambda^2}: 
(\Lambda^2,\bold C \omega_X')\rightarrow (\Lambda^2, \bold C \omega_X')$, 
where $\Lambda^2=\Lambda_{\text{K3}}$.  

Recall that the reflections $r_{C_1}, \cdots, r_{C_m}$ for 
smooth $(-2)$ curves \break $C_1,\cdots,C_m$ come from 
the twisting functors $\Phi^{\Cal O_X(C_1)}_X \circ T_{C_1}, \cdots, 
\Phi^{\Cal O_X(C_m)}_X \circ T_{C_m}$.  
These reflections act as identity on $\Lambda^0, \Lambda^4$ and act on 
$\Lambda^2$.  The composition 
$$
l:=r_{C_1}\circ \cdots \circ r_{C_m} \circ k : 
(\tilde \Lambda, \bold C \omega_X') \rightarrow 
(\tilde \Lambda, \bold C \omega_X') 
$$
satisfies  
$l((0,0,1))=(0,0,1), l((1,0,0))=(1,0,0)$ and $l: 
(\Lambda^2, \bold C \omega_X') \rightarrow 
(\Lambda^2, \bold C \omega_X')$. 
Now define the positive cone $\Cal P^+$ to be one of the two connected 
components of 
$\{ x \in \tau_X(NS(X)) \otimes \bold R \vert (x,x) > 0 \}$ 
which contains the ample class. Using the Nakai-Moishezon criterion, 
the Hodge index Theorem, and the fact that any irreducible curve $C \subset 
X$ satisfies $C^2 \geq 0$ unless $C\cong \bold P^1$, we may describe the 
ample cone $A(X)$ by 
$$
A(X)=\{ x \in \Cal P^+ \;\vert\; (x,C)>0, \forall C \cong \bold P^1 \}. 
$$
This means $A(X)$ is the fundamental domain of the reflection group 
$\langle r_{C} \vert C \cong \bold P^1, C \subset X \rangle$ acting 
on $\Cal P^+$. 
Therefore composing suitable reflections, 
and $\iota_2$ if necessary\footnote{ We thank B. Szendr\"oi for 
pointing this out to us.}, 
we can assume that $l$ preserves the ample $A(X)$. 
(This is the only place where we may need the Hodge isometry $\iota_2$.) 
Therefore there is 
an automorphism $\varphi \in \text{Aut}(X)$ such that $l=\varphi^* (=
f^{\Cal O_{\Gamma(\varphi)}}_X)$ by the global Torelli theorem. 
Here we note that $\varphi^*$ is identity on $\Lambda^0$ and $
\Lambda^4$. To summarize, we have obtained an identity 
$$
r_{C_1}\circ \cdots \circ r_{C_m} \circ (\iota_2)^s\circ 
f_X^{(-B)}\circ 
(F(f_X^{\Cal L_1}, f_X^{\Cal L_2}, f^{\Cal I_\Delta}_X)\circ g)^{-1}
\circ f^{\Cal E}_{X\rightarrow X}= \varphi^*
$$
as an element of $O_{Hodge}((\tilde \Lambda,\bold C \omega_X'))=
O_{Hodge}(\tilde H(X,\bold Z))$. In the above equation all, except possibly 
$g$, are in the group generated by the image of the homomorphism 
$$
ch: \text{Auteq} D(X) \rightarrow 
O_{Hodge}(\tilde H(X,\bold Z))=O_{Hodge}(\tilde \Lambda,\bold C \omega_X'), 
$$
and $\iota_2$. 
Therefore we conclude that $g \in \langle \text{Im}(ch), \iota_2 \rangle$ 
as well.  \qed

\head
{\S 3 Symplectic mapping class group and Proof of Theorem \MainThII }
\endhead

\global\secno=3  
\global\propno=1 
\global\eqnum=1 

In this section we prove the surjectivity of the map 
$Symp(Y,\kappa_Y)$  $\rightarrow O^+(T(Y))^*$ which implies Theorem \MainThII. 
(Note that $Symp^0(Y,\kappa_Y)$ acts on $H^2(Y,\bold Z)$ trivially.)

Let $Y$ be a K3 surface and $A(Y) \subset NS(Y)\otimes \bold R$ be the ample 
cone. As described in Definition {\sympY}, for $\kappa_Y \in A(Y)$ 
we consider a symplectic structure $(Y,\kappa_Y)$ and the group of symplectic 
diffeomorphisms $\text{Symp}(Y,\kappa_Y)$. We assume this symplectic structure 
$(Y,\kappa_Y)$ is generic. Under this assumption we consider the induced action 
$g^*$ on $H^2(Y,\bold Z)$ of $g \in \text{Symp}(Y,\kappa_Y)$. 

\Proclaim{Lemma}{\LemmaI} Let $(Y,\kappa_Y)$ be a generic symplectic structure 
in the sense of Definition {\sympY}, 4)  
and $g\in \text{Symp}(Y,\kappa_Y)$. Then we have
\item{1)} $g^*(NS(Y))=NS(Y)$, in fact $g^*\vert_{NS(Y)}=\text{id}_{NS(Y)}$,
\item{2)} $g^*(T(Y))=T(Y)$.
\endproclaim

\demo{Proof} Let us first recall that if $\kappa_Y$ is generic, then 
$NS(Y)$ is characterized as the minimal primitive sublattice of 
$H^2(Y,\bold Z)$ which 
contains $\kappa_Y$ after tensoring with $\bold R$. (See the definition 
and Remark given after Definition {\sympY}.) This implies  
$g^*(NS(Y))=NS(Y)$, since $g^*(\kappa_Y)=\kappa_Y$ and $g^*(H^2(Y,\bold Z)) = 
H^2(Y,\bold Z)$.  (Note that a diffeomorphism $g \in Diff(Y)$ does not 
preserve the lattice $NS(Y)$ in general.) The induced action $g^*$ of $g \in 
Diff(Y)$ on $H^2(Y,\bold Z)$ is an isometry. Therefore we may conclude 
$g^*(T(Y))=T(Y)$ as well from $g^*(NS(Y))=NS(Y)$ and $T(Y)=NS(Y)^\perp$ 
in $H^2(Y,\bold Z)$.  Moreover, one can show 
$g^*\vert_{NS(Y)}=\text{id}_{NS(Y)}$ as follows: Let $x \in NS(Y)$. Then 
$(\kappa_Y,g^*x)=(g^*\kappa_Y,g^*x)=(\kappa_Y,x)$, i.e. 
$(\kappa_Y,g^*x-x)=0$ and $g^*x-x \in NS(Y)$. If $g^*x-x\not=0$, then 
one has $\kappa_Y \in (g^*x-x)^\perp \otimes \bold R$, where $(g^*x-x)^\perp$ 
is the orthogonal lattice in $NS(Y)$. This is a contradiction to the 
minimality of $NS(Y)$ for generic $\kappa_Y$. This proves $g^*x=x$ 
for all $x \in NS(Y)$.  \qed 
\enddemo

Recall the definition of the kernel subgroup $O(T(Y))^*$ and its 
index two subgroup $O^+(T(Y))^*=O(T(Y))^* \cap O^+(T(Y))$. We can now 
prove our Theorem {\MainThII} using the Theorem {\Donaldson}. 

\vskip0.2cm

\noindent
{\it Proof of Theorem {\MainThII}:} Let $g \in \text{Symp}(Y,\kappa_Y)$. 
Consider the induced action on the discriminant groups 
$A_{NS(Y)}=NS(Y)^*/NS(Y)$ and $A_{T(Y)}=T(Y)^*/T(Y)$. 
Then by Lemma {\LemmaI} we see that $NS(Y)$ and 
$T(Y)$ are stable under the induced action $g^*$, and $g^*\vert_{A_{NS(Y)}} 
=\text{id}_{A_{NS(Y)}}$. This implies $g^*\vert_{A_{T(Y)}}= 
\text{id}_{A_{T(Y)}}$, and $g^*\vert_{T(Y)} \in O(T(Y))^*$. 
Since $g \in \text{Symp}(Y,\kappa_Y) \subset Diff(Y)$, we have 
$g^* \in O^+(H^2(Y,\bold Z))$ by Theorem {\Donaldson} and 
$g^*(\kappa_Y)=\kappa_Y$  by the assumption. 
Putting all together, we conclude 
that $g^*\vert_{T(Y)} \in O^+(T(Y))$ and further 
$g^*\vert_{T(Y)} \in O^+(T(Y))^*$. Other statements in Theorem {\MainThII}, 
1) are included in Lemma {\LemmaI}. 

To show the surjectivity of the map $\text{Symp}(Y,\kappa_Y) \rightarrow 
O^+(T(Y))^*$, let us take $\varphi \in O^+(T(Y))^*$. Then $(\text{id}, 
\varphi) \in O(NS(Y)) \times O(T(Y))$ extends to $\tilde \varphi \in 
O(H^2(Y,\bold Z))$ since $\text{id}|_{A_{NS(Y)}}=
\text{id}_{A_{NS(Y)}}$ and   $\varphi|_{A_{T(Y)}} = \text{id}_{A_{T(Y)}}$. 
This extension $\tilde \varphi$ is an element in $O^+(H^2(Y,\bold Z))$ 
since $\varphi \in O^+(T(Y))$ and $\text{id}(\kappa_Y)=\kappa_Y$. Then 
Theorem {\Donaldson} ensures the existence 
$g \in Diff(Y)$ such that $g^* = \tilde \varphi$. By construction of 
$\tilde \varphi$, we have $g \in \text{Symp}(Y,\kappa_Y)$ and 
$g^* \vert_{T(Y)}=\varphi$. This shows the surjectivity. \qed

\head
{\S 4  Monodromy group, FM partners, and Proof of Theorem {\MainThIII}} 
\endhead

\global\secno=4  
\global\propno=1 
\global\eqnum=1 

\vskip0.2cm

Here we present Proof of Lemma {\genericLemma} and Proof of 
Theorem {\MainThIII}. 

\vskip0.3cm

\noindent
{\it Proof of  Lemma {\genericLemma}:} Recall that 
for generic $(\vX_1,\tau_1)$ and 
$(\vX_2,\tau_2)$ we have $\tau_i(T(\vX_1))$ $=$ $U\oplus M_n$ for $i=1,2$. 
If  $\vX_1 \cong \vX_2$, then there exists a Hodge isometry 
$\varphi^*: H^2(\vX_2,\bold Z) \rightarrow H^2(\vX_1,\bold Z)$.  
Since the Hodge isometry maps the transcendental lattice to 
the transcendental lattice, i.e. $\varphi^*: T(\vX_2) \rightarrow 
T(\vX_1)$, the restriction of $\tau_2\circ (\varphi^*)^{-1}\circ \tau_1^{-1}$ 
to $U \oplus M_n$ gives a desired isometry $g$ of $U \oplus M_n$. 

Conversely, assume there exists $g \in O(U\oplus M_n)$ with the property  
$g(\bold C \tau_1(\omega_{\vX_1}))=\bold C \tau_2(\omega_{\vX_2})$. 
Since the primitive embedding $U \oplus M_n \hookrightarrow 
\Lambda_{\text{K3}}$ is unique up to isomorphism [Ni,(1.14.4)], 
there is an isometry $\bar g: 
\Lambda_{\text{K3}}\rightarrow \Lambda_{\text{K3}}$ such that $\bar g 
\vert_{U \oplus M_n} = g$. Then $\tau_2^{-1}\circ \bar g \circ \tau_1$ 
is a Hodge isometry from $H^2(\vX_1,\bold Z)$ to $H^2(\vX_2,\bold Z)$. 
Therefore $\vX_1 \cong \vX_2$ by the Torelli Theorem. \qed 

\vskip0.3cm
We can now proceed to Proof of Theorem {\MainThIII}. However we 
state one general result whose special case will be used in the proof. 

\Proclaim{Lemma}{\Niklin}  
Let $M$ be a lattice of signature $(1,t)$ and 
$U$ be the hyperbolic lattice. Consider the lattice $N=U\oplus M$. Denote the 
discriminant by $A_N=N^*/N=M^*/M$ and consider the orthogonal group 
$O(A_N)$ with respect to the natural induced form on $A_N$. 
Then 
$$
[O(N):O(N)^*]=|O(A_N)|, 
$$
where $O(N)^*=\text{Ker}( O(N) \rightarrow O(A_N) )$.
\endproclaim
\demo{Proof} Since $N=U\oplus M$, we have $l(A_{N})=l(A_M)$ and 
$\text{rk} N \geq l(A_N) + 2$. 
Here $l(A)$ denotes the minimal number of generators of a finite 
abelian group $A$. 
Therefore by [Ni,Theorem 1.14.2], 
the natural map $O(N) \rightarrow O(A_N)$ is surjective.  
Then the result follows immediately 
by the definition of the kernel $O(N)^*$. \qed 
\enddemo

\noindent
{\it Proof of Theorem {\MainThIII}: } Since the case $n=1$ is easy, we 
consider only the case $n\geq2$. 
For the lattice $M_n = \langle 2n \rangle$ the group $O(A_{U\oplus M_n})$ 
is easily described by the Chinese remainder theorem as 
$$
O(A_{U\oplus M_n})=O(M_n^*/M_n)=O(\bold Z / 2n) \cong (\bold Z/2)^{p(n)},
$$
where $p(n)$ is the number of the primes dividing $n$ (cf. 
[Sc, lemma 3.6.1] and the proof there). 
Using the Lemma {\Niklin} we have 
$$
[O(U\oplus M_n): O(U\oplus M_n)^*]=|O(A_{U\oplus M_n})|=2^{p(n)} .
$$
Since $O(U\oplus M_n)^* \rightarrow 
O(U\oplus M_n)/\{ \pm \text{id} \}$ 
is injective (Lemma {\compo}), we have 
$$
[O(U\oplus M_n)/\{ \pm \text{id} \} : O(U\oplus M_n)^*]
= 2^{p(n)-1} .
$$
By definition we have $[O(U\oplus M_n)^*: O^+(U\oplus M_n)^*]\leq 2$,  
and we find an element $(-\text{id}_U,\text{id}_{M_n})$ in $O(U\oplus M_n)^* 
\setminus O^+(U\oplus M_n)^*$, thus we can conclude 
$[O(U\oplus M_n)^*:O^+(U\oplus M_n)^*]=2$. 
We have $[O(U\oplus M_n)/\{ \pm \text{id} \}: 
O^+(U\oplus M_n)/\{ \pm\text{id} \}]=2$ as well. 
Therefore we finally conclude 
$$
[O^+(U\oplus M_n)/\{ \pm \text{id} \} : O^+(U\oplus M_n)^*]
=2^{p(n)-1},
$$
i.e. $[\Cal{M}_n(\vX): \Cal{MS}_n(\vX) ]= 2^{p(n)-1}$ \qed

\head
{\S 5 Mirror family of a K3 surface with $\text{deg}(X)=12$ } 
\endhead

\global\secno=5  
\global\propno=1 
\global\eqnum=1 
\vskip0.2cm

\noindent
{\it 1) The mirror family:}  
Consider the following one-parameter family of a surface given by 
three quadrics and a hyperplane in $\bold P^6$;
$$
Q(\psi): 
\cases
U_1+U_2+U_3+U_4+U_5+U_6- \psi U_0 =0 &  
( \psi \in \bold C \setminus\{ \pm2, \pm 6\} ) \cr
U_1U_2=U_0^2, \;\;\; U_3U_4=U_0^2, \;\;\; U_5U_6 = U_0^2, &  \cr 
\endcases
$$
where $U_0,\cdots, U_6$ are homogeneous coordinates of $\bold P^6$. 
$Q(\psi)$ has 12 double points at $U_0=0$, which may be blown up to 
$$
\tilde Q(\psi): 
\cases
U_1+U_2+U_3+U_4+U_5+U_6- \psi U_0 =0 &  ( \psi\not=\pm2, \pm 6) \cr
U_1V_1=V_2U_0 ,\;\;\; U_3V_3=V_4U_0, \;\;\; U_5V_5=V_6U_0 & \cr
U_2V_2=V_1U_0, \;\;\; U_4V_4=V_3U_0, \;\;\; U_6V_6=V_5U_0, & \cr
\endcases 
$$
where $([U_0,\cdots,U_6],[V_1,V_2],[V_3,V_4],[V_5,V_6]) 
\in \bold P^6 \times \bold P^1 \times \bold P^1 \times \bold P^1$. 
The natural map $\tilde Q(\psi) \rightarrow Q(\psi)$ is a crepant 
resolution and we see that $\tilde Q(\psi)$ is a smooth 
K3 surface. Also we may verify that the singularities of 
$\tilde Q(\psi)$ at $\psi=\pm2, \pm 6$ are ordinary double points.

\Proclaim{Theorem}{\PeSt}{\bf ([PS, Theorem 1])} For generic 
$\psi \in \bold C\setminus\{\pm2,\pm6\}$, one has 
$$
\text{NS}(\tilde Q(\psi))\cong E_8(-1)^{\oplus 2} \oplus U \oplus 
\langle -12 \rangle .
$$
\endproclaim
\demo{Sketch of Proof} Blowing up the 12 double points introduces 
the following 12 lines; i.e. a line 
$$
l_{0,+1,+1}:=
\{ ([0,0,0,1,0,-1,0],[V_1,V_2],[0,1],[0,1]) \vert (V_1,V_2) \in \bold P^1 \}
$$
and its obvious permutations of the coordinates, which appear 
12 middle points of the cube in Fig.1. As we see in Fig.1, 
there are 8 lines of different type. These are represented by lines 
(corresponding to the corners of the cube in Fig.1), e.g.
$$
l_{+1,+1,+1}:= 
\{ ([0,U_1,0,U_3,0,U_5,0],[0,1],[0,1],[0,1]) \vert 
U_1+U_3+U_5=0 \}.
$$
These two types of lines are independent of $\psi$ and altogether 
generate a sublattice of rank 17 in $NS(\tilde Q(\psi))$. 
What is interesting, and even crucial for the following analysis, is that 
we have 12 more lines which {\it dependents} on $\psi$, and 
which are typically represented by 
$$
\aligned
m_{1,+1,+1}(\beta):=
\{ ([U_0,\beta U_0, \beta^{-1} U_0, U_3, U_4, -U_3, &-U_4],[1,\beta],
 [V_3,V_4],[V_3,-V_4]) \; \vert \; \cr 
& U_3 V_3 = V_4 U_0, U_4 V_4=V_3 U_0 \} 
\endaligned
$$ 
with $\beta$ satisfying $\beta+\beta^{-1}=\psi$. Suitable permutations 
of the coordinates generate 12 lines of this type. These lines 
contribute two more linearly independent classes in $NS(\tilde Q(\psi))$, 
and in total we see $\text{rk} \tilde Q(\psi)\geq19$.  
On the other hand, 19 is the maximal Picard number for a non-constant 
family of K3 surfaces (see e.g. [Og1, Main Theorem]). 
Therefore $\text{rk} \tilde Q(\psi)=19$ for generic $\psi \in \bold C\setminus
\{ \pm 2, \pm 6 \}$. 
Working out the intersection numbers in 
detail, we obtain the desired result. \qed  
\enddemo

\vskip0.5cm

\centerline{\epsfxsize 2.5truein\epsfbox{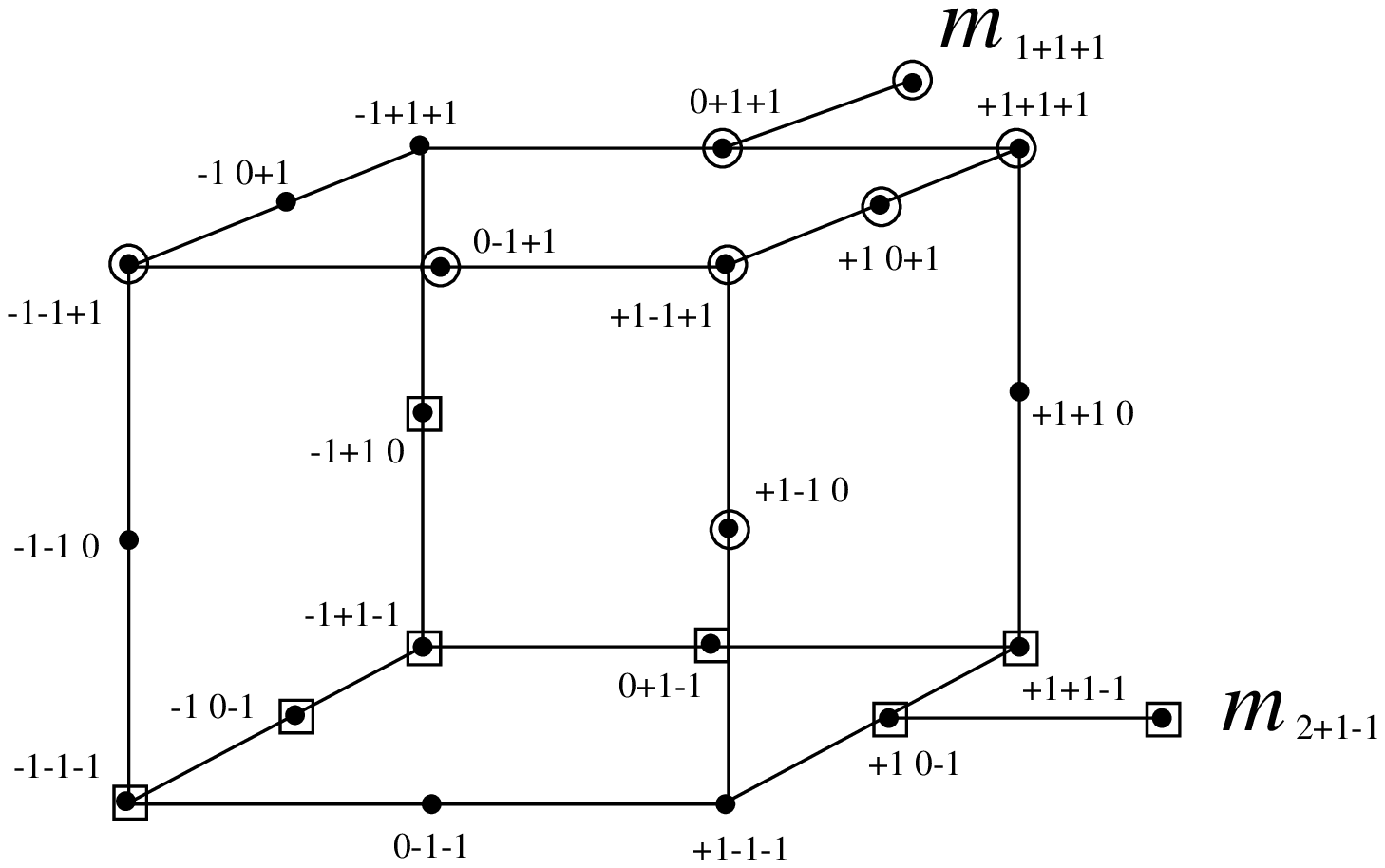}}

\botcaption{Fig.1. Cube $[-1,1]^3$ representing the lines}
The notation $+1+1-1$, for example, represents the line $l_{+1,+1,-1}$. 
Two additional lines, $m_{1,+1,+1}, m_{2,+1,-1}$, are also written to 
show the generators of two $E_8(-1)$  lattices. 
\endcaption

\vskip0.5cm

\noindent
{\bf Remark.} $\tilde Q(\psi)$ may be written in the affine chart $U_0\not=0$ 
by a single equation,
$$
f(X,Y,Z)=X+{1\over X}+Y+{1\over Y}+Z+{1\over Z} - \psi,
\Eqno{\defeq}
$$
where $X=U_1/U_0,Y=U_3/U_0,Z=U_5/U_0$. The Newton polytope of this 
defining equation is simply given by the convex hull 
$\text{Conv.}\big( \{ (\pm1,0,0), (0,\pm1,0),(0,0,\pm1)\} \big)$. 
Then the polar dual 
of this polytope is given by $[-1,1]^3$, which is the cube depicted in Fig.1. 
In fact, the polytope $[-1,1]^3$ is so-called {\it reflexive polytope} 
considered 
by Batyrev [Ba1] to explain the mirror symmetry of Calabi-Yau hypersurfaces 
in toric varieties. If we follow his construction, we obtain a 
defining equation 
$$
f(X,Y,Z)=\psi_1(X+{1\over X}) + \psi_2 (Y +{1\over Y}) + \psi_3 (Z+{1\over Z})
+ 1
$$
for a generic hypersurface in a toric (Fano) variety whose crepant resolution 
may be read from the toric diagram (the cube) in Fig.1. When we fix a 
crepant resolution of the ambient toric variety, the generic 
hypersurface defines a family of smooth K3 surfaces with the N\'eron-Severi 
lattice of rank 17, which is generated by the lines appeared 
in the vertices of the cube. Our K3 surface 
$\tilde Q(\psi)$ is a specialization of this generic family to $\psi_1=\psi_2=
\psi_3=-{1\over \psi}$ (see [LY2]).  \qed

\Proclaim{ Proposition}{\mirrorMn} $\{ \tilde Q(\psi) \}_{  
\psi \in \bold C \setminus \{ \pm2,\pm6\}}$ 
is a (covering of the) family mirror to $\langle 12\rangle$-polarizable 
K3 surface. 
\endproclaim
\demo{Proof} For the lattice $M_n=\langle 2n \rangle$, the orthogonal 
lattice $M_n^\perp$ in the K3 lattice $\Lambda_{\text{K3}}$ has the form 
$M_n^\perp=U\oplus \vM_n$ with 
$$
\vM_n = \langle -2n \rangle \oplus U \oplus E_8(-1)^{\oplus 2}.
$$
Then combining this with Theorem {\PeSt}, we conclude the 
statement for $n=6$. \qed
\enddemo 

\noindent
{\bf Remark.} Note the isomorphism $\tilde Q(\psi)\cong \tilde Q(-\psi)$. 
Indeed, in the formula (\defeq) one sees a birational map $\tilde Q(\psi) 
\rightarrow \tilde Q(-\psi): (X,Y,Z) \rightarrow (-X,-Y,-Z)$. 
Since any birational map between K3 surfaces is biregular (isomorphism), 
the claim follows.  \qed

\vskip0.5cm
\noindent
{\it 2)  Period integrals, PF equation and monodromy: }
For $|\psi|>6$ consider the integral;
$$
\Pi(\psi)={1\over (2\pi i)^3} \int_{|X|=|Y|=|Z|=1} {-\psi \over f(X,Y,Z)} 
{dX \over X} 
{dY \over Y} 
{dZ \over Z} 
$$
for the defining equation $f(X,Y,Z)$ (\defeq) of $\tilde Q(\psi)$ in the 
affine chart $U_0\not=0$. This integral represents a period integral 
for the family $\tilde Q(\psi)$ (see [PS],[Ba2]). 

\Proclaim{Proposition}{\PFeq} {\bf ([Ba2], [LY2])} 
When $|\psi|>6$, the integral can be evaluated to be a convergent power series 
$$
\Pi(x)=\sum_{k,l,m \geq 0} c(k,l,m) x^{k+l+m}:=
\sum_{k,l,m \geq 0} {(2(k+l+m))! \over (k!)^2(l!)^2 (m!)^2 } 
x^{k+l+m} \quad (x:={1\over \psi^2}),
$$
and satisfies the following differential equation:
$$
\{\theta^3+36x^2(\theta+1)(2\theta+1)(2\theta+3) 
-2x(2\theta+1)(10\theta^2+10\theta+3)\} \Pi(x) =0,
\Eqno{\PF}
$$
where $\theta=x {d \; \over dx}$.  \qed
\endproclaim

The differential equation is defined on $\bold P^1$ 
and has regular singularities at $x=0,{1\over 36}, 
{1\over 4}, \infty$. Thus the affine space $\bold C$ 
is naturally compactified to $\bold P^1$ considering 
the boundary point $\psi=\infty \; (x={1\over \psi^2}=0)$. 
(Note that $\pm \psi$ are naturally identified in accord with  
the isomorphism $\tilde Q(\psi)\cong \tilde Q(-\psi)$.) 
The boundary point $x=0$ 
is the so-called {\it maximally unipotent monodromy point} [Mo1], and 
plays important roles in the applications of mirror symmetry. As we see 
in (\PF), about this boundary point all indices of the local solutions 
degenerate to zero, and we may apply the Frobenius method to 
generate all solutions (see, e.g. [HLY] and references therein for details). 
Now we define a ratio of two solutions
$$
t:={1 \over 2\pi i} {\partial \; \over \partial \rho} 
\log \Pi(x,\rho)\vert_{\rho=0}
={1\over 2\pi i} {1\over \Pi(x)} {\partial \; \over \partial \rho} 
\Pi(x,\rho)\vert_{\rho=0},
\Eqno{\mirrormap}
$$
where $\Pi(x,\rho):=\sum_{k,l,m\geq0} c(k+\rho,l+\rho,m+\rho)x^{k+l+m+3\rho}$ 
is the formal extension of the series to apply the Frobenius method. 
The ratio ({\mirrormap}) may be regarded as a projective coordinate of 
the period point in $\Omega(U\oplus M_n)$ (cf. [Do, Section 7]). 
The inverse relation to ({\mirrormap}), i.e. $x=x(t)$ is called the 
{\it mirror map}. In [LY1], nice modular properties related to the genus 
zero modular groups have been observed for many examples of 
the mirror family of $M_n$-polarized K3 surfaces. 
In our case of $\tilde Q(\psi)$, it may be summarized as follows:

\Proclaim{Proposition}{\LYmodular} {\bf ([LY2])} The mirror map 
$x=x(t)$  is the Hauptmodul of the genus zero modular group 
$\Gamma_0(6)_+$ and satisfies the following Schwarzian equation: 
$$
\{ t,x \} = {1-52x+1500x^2-6048x^3+15552 x^4 \over 2(1-36x)^2(1-4x)^2 x^2}, 
\qquad 
$$
where $\{ t, x \}={t''' \over t'}-{3\over2}({t'' \over t'})^2$ with 
$':={d \; \over dx}$. \qed
\endproclaim

\eject

\noindent
{\bf Remark.}  1) We follow the notation in [CN] for the extension of the 
group $\Gamma_0(n)$ by the involutory normalizers, the Fricke involution and 
the Atkin-Lehner involutions. The Fricke involution is a normalizer 
in $PSL(2,\bold R)$ of the modular subgroup, 
$$
\Gamma_0(n)=\big\{ \left( \matrix a & b \cr c & d \cr \endmatrix \right) \in 
PSL(2,\bold Z) \; \vert \; c\equiv 0 \; \text{mod} n \; \big\}, 
$$
and it is represented by a coset 
$
W_n:= \left( 
\matrix 0 & -{1\over \sqrt{6}} \cr \sqrt{6} & 0 \cr \endmatrix \right)
\Gamma_0(n) .
$
The Atkin-Lehner involutions $W_{r}$ are generalizations of the 
Fricke involution $(r=n)$ where  $r \geq 1$ is a divisor of $n$ such that 
$r$ and $n/r$ are coprime. ($W_1=\Gamma_0(n)$.) Then the notation 
$\Gamma_0(n)_+$ is used to represent the group obtained from $\Gamma_0(n)$ 
by adjoining all possible Atkin-Lehner involutions, while the 
$\Gamma_0(n)_{+n}$ represents the group joined the Fricke involution only. 
In our case we may work out explicit generators of the groups as follows;
$$
\aligned
& 
\Gamma_0(6)_+=\big\langle  
\left( \matrix 1 & 1 \cr 0 & 1 \cr \endmatrix \right), 
\left( \matrix 0 & -{1\over \sqrt{6}} \cr \sqrt{6} & 0 \cr \endmatrix \right), 
\left( \matrix \sqrt{3} & {1\over \sqrt{3}} \cr 2\sqrt{3} & \sqrt{3} 
       \cr \endmatrix \right) 
\big\rangle , 
\cr
& 
\Gamma_0(6)_{+6}=\big\langle  
\left( \matrix 1 & 1 \cr 0 & 1 \cr \endmatrix \right), 
\left( \matrix 0 & -{1\over \sqrt{6}} \cr \sqrt{6} & 0 \cr \endmatrix \right), 
\left( \matrix 5 & 2 \cr 12 & 5 \cr \endmatrix \right) 
\big\rangle .
\cr
\endaligned
\Eqno{\Gsubs}
$$

\noindent
2) In [BP] a different family of K3 surfaces and therefore a different 
Picard-Fuchs equation is studied. The family in [BP] is a different 
covering of the mirror family of $\langle 12 \rangle$-polarizable 
K3 surfaces.  The definitions of the Period integral 
and the mirror map $x=x(t)$ are parallel to our family $\tilde Q(\psi)$. 
In this case, we obtain for the mirror map $x=x(t)$ the Hauptmodul  
of the genus zero modular group $\Gamma_0(6)_{+6}$, which is an index two 
subgroup of $\Gamma_0(6)_+$. The difference of the modular group comes 
from the difference of the covering. In fact, the precise relation 
between these two parametrizations may be found in [PS, Remark4]. 

\noindent
3) The Schwartzian equation in Theorem {\LYmodular} has four singularities 
at $0,{1\over36},{1\over4},\infty$. To see these singularities in one affine 
chart, it is convenient to apply a fractional linear transformation $z={48 x 
\over 12 x +1}$ which sends the four singular points 
$(0,{1\over36},{1\over4},\infty)$ to $(a_0,a_1,a_2,a_3)=(0,1,3,4)$. 
Then the Schwartzian equation may be written in the standard form 
(see [Hi, Chapter 10] for example): 
$$
\{ t ,z \} = \sum_{i=0}^3 \big( {1\over2}{1-\alpha_i^2 \over (z-a_i)^2} 
+ {\beta_i \over z-a_i } \big), 
$$
with $(\alpha_0,\alpha_1,\alpha_2,\alpha_3)=(0,{1\over2},{1\over2},{1\over2})$ 
and $(\beta_0,\beta_1,\beta_2,\beta_3)=({13\over24},-{3\over16},{1\over48}, 
-{3 \over8})$. The $\alpha_i$'s determine the local form of the mapping 
$t(z)=c_0 (z-a_i)^{\alpha_i}+\cdots $ near the singularities. 
From this we can determine the image of the complex $x$ plane to the upper 
half plane parametrized by $t$. In fact, since the group $\Gamma_0(6)_+$ is 
a genus zero group, the image coincides with a fundamental domain of 
the group. In Fig. 2, we have depicted the fundamental domain of 
$\Gamma_0(6)_+$ and also that of $\Gamma_0(6)_{+6}$.  \qed

\vskip0.5cm

\centerline{\epsfxsize 4truein\epsfbox{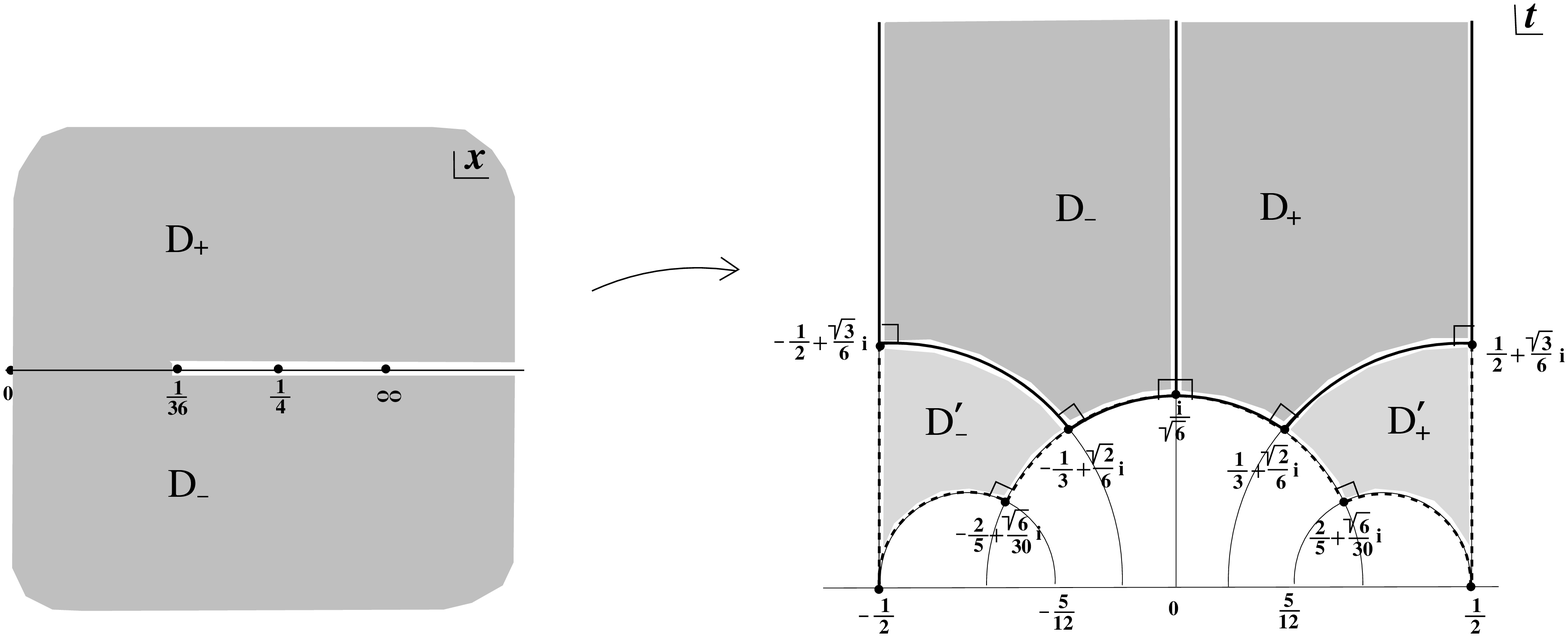}}

\botcaption{Fig.2.  Fundamental domain of $\Gamma_0(6)_+$}
The fundamental domain of $\Gamma_0(6)_+$ is shown as an image of 
a complex $x$-plane by $t=t(x)$. See also Table 1 for the correspondence  
of the (regular) singular points and the points of non-trivial stabilizer 
on the $t$-plane.  The union $D_+\cup D_- \cup D_+' \cup D_-'$ represents 
the fundamental domain of $\Gamma_0(6)_{+6}$. 
\endcaption

\vskip0.5cm

$$
\def\vspace#1{  height #1 & \omit && \omit && \omit && \omit &  \cr }
\vbox{\offinterlineskip
\hrule
\halign{ \strut
\vrule#&  $\;$ \hfil $#$ \hfil 
&\vrule#&  $\;$ \hfil $#$ \hfil  
&\vrule#&  $\;$ \hfil $#$ \hfil  
&\vrule#&  $\;$ \hfil $#$ \hfill  
&\vrule# 
\cr 
&  x && t &&  \text{stabilizer}  &&  \text{FM transform}  &\cr
\noalign{\hrule} 
\noalign{\hrule} 
&  0 
&&  i \infty 
&&  T=\left( \matrix 1 & 1 \cr 0 & 1 \cr \endmatrix \right)  
&&  \Phi_0=\Phi_X^{\pi_2^* \Cal O(1)} 
&\cr
&  {1\over 36} 
&& {i\over \sqrt{6}} 
&&   S_1= 
     \left( \matrix 0 & {-1\over\sqrt{6}} \cr 
            \sqrt{6} & 0 \cr \endmatrix \right)
&&   \Phi_1=T_{\Cal O_X}  
&\cr
&   {1\over 4} 
&&  {1 \over3}+{\sqrt{2}\over 6}i  
&&   S_2 =  
     \left( \matrix -\sqrt{2} & {1\over\sqrt{2}} \cr 
                  -3\sqrt{2} & \sqrt{2} \cr \endmatrix \right)
&&  \hskip0.5cm  -  
&\cr
&   \infty 
&&  {1\over2}+{\sqrt{3} \over 6}i  
&&   S_2S_1 T^{-1} 
&&  \Phi_0^{-1}\circ \Phi_1 \circ ``(-)'' 
&\cr
}
\hrule}
$$

\botcaption{ Table 1. Monodromy matrices and FM transforms} 
For the regular singular points on the $x$-plane, we list the 
corresponding monodromy matrices as an element in $\Gamma_0(6)_+$ 
and FM transforms. The blank (-) in the third (and the fourth) line 
indicates that this does not come from $\text{Auteq} D(X)$ 
(see Proposition 5.8 and Remark after that). 
\endcaption

\vskip0.5cm

The relation between the modular group and the monodromy group of 
the period integrals is described by the well-known relation;  
$$
PSL(2,\bold R) \cong SO^+(2,1; \bold R).   
\Eqno{\PSLisoPSO}
$$
Let us describe this isomorphism explicitly for the Minkowski space 
$\bold R^{2,1}=(U\oplus M_n)\otimes \bold R$ fixing a basis $e, v, f$ for 
the lattice $U\oplus M_n$ ( precisely $U=\bold Z e \oplus \bold Z f$ and 
$M_n=\bold Z v$ with $(v,v)=2n$). Then 
the orthogonal group is defined by
$$
SO(2,1;\bold R):=\{ g \; \vert \; \,^tg \Sigma g  = \Sigma, \text{det}(g)=1 \} 
\;, \; 
\Sigma=
\left(
\matrix 0 & 0 & -1 \cr 
0 & 2n & 0 \cr
-1 & 0 & 0 \cr
\endmatrix\right) .
$$
$SO^+(2,1; \bold R)$ is an index two subgroup consisting of the elements 
which preserve orientations of the positive two planes.  We can write 
the isomorphism ({\PSLisoPSO}) (anti-homomorphism) explicitly; 
$$
R: \left( 
\matrix 
a & b \cr 
c & d \cr 
\endmatrix \right) 
\mapsto \left( 
\matrix 
a^2 & 2 a c & {c^2 \over n} \cr
ab & ad+bc & {cd \over n} \cr
n b^2 & 2n b d & d^2  \cr 
\endmatrix \right). 
\Eqno{\Rmap}
$$
We may also consider the natural map $O(U\oplus M_n) \rightarrow 
O(2,1;\bold R)$ and define a homomorphism 
$F: O^+(U\oplus M_n) \rightarrow SO^+(2,1; \bold R)$ by 
$F(g)=\text{det}(g) g$. Then we have  
$\text{Ker} F=\{ \pm \text{id}_{U\oplus M_n} \} $ and thus 
the induced map 
$$
\bar F: {\Cal M}_n(\vX)=O^+(U\oplus M_n)/\{ \pm \text{id}_{U\oplus M_n} \} 
\rightarrow SO^+(2,1;\bold R)\cong PSL(2,\bold R)
$$
is injective. For a subgroup $G \subset O^+(U\oplus M_n)$, 
consider the following composition map: 
$$
q: G \hookrightarrow O^+(U\oplus M_n) \rightarrow 
O^+(U\oplus M_n)/\{ \pm \text{id}_{U\oplus M_n} \}. 
$$
Then, up to the kernel $\text{Ker} (q)$, we can identify $G$ with the image 
$R^{-1}\circ \bar F\circ q (G)$ in $PSL(2,\bold R)$. 
For example, using the (injective) map ({\compOp}) 
we may regard the group $O^+(U\oplus M_n)^*$ $(n\geq 2)$ as a subgroup of 
$PSL(2,\bold R)$. 

The following group isomorphisms are derived in purely arithmetic 
manner in [Do], and will be utilized in the following. 

\Proclaim{Theorem}{\Do} {\bf ( [Do, Theorem (7.1), Remark (7.2)])}  We have: 
\item{1)} $O^+(U\oplus M_n)^* \cong \Gamma_0(n)_{+n} \;\; (n\geq 2)$, 
$O^+(U\oplus M_1)^*/\{ \pm \text{id}_{U\oplus M_1} \} \cong PSL(2, 
\bold Z)$ 
\item{2)} $O^+(U\oplus M_n)/\{ \pm \text{id}_{U\oplus M_n} \}  
\cong \Gamma_0(n)_+$ 
\qed
\endproclaim

\noindent
{\bf Remark.} When $\text{Ker}(q)$ is non-trivial, we 
encounter the sign ambiguity to determine $G \subset O^+(U\oplus M_n)$ 
from its image in $PSL(2,\bold R)$.  However 
the different choice of the sign can be understood as coming from 
the Hodge isometry associated to the shift functor $[1]$ or 
its counterpart in $\text{Auteq} D Fuk(\vX)$. Because of this,  
we do not have to pay much attention to this sign problem. 
This is a reason we define the monodromy group $\Cal M_n(\vX)$ 
by the quotient $O^+(U\oplus M_n)/\{ \pm \text{id}_{U\oplus M_n} \}$ 
instead of the conventional definition of the monodromy group 
of the period integrals.  \qed

\vskip0.5cm
\noindent
{\it 3) Monodromy $S_2$ and FM partner: } 
For each stabilizer $T,S_1,S_2 \in \Gamma_0(6)_+$ defined in Table 1, 
consider the matrices $R(T),R(S_1), R(S_2) \in SO^+(2,1;\bold R)$. We can 
verify these matrices are in the image of the map $O^+(U\oplus M_n)/\{ \pm 
\text{id}_{U\oplus M_n} \} \rightarrow SO^+(2,1;\bold R)$. 
Though the preimages $\bar T, \bar S_1, \bar S_2 \in O^+(U\oplus M_n)$ 
of $R(T), R(S_1), R(S_2)$ are determined only up to signs as matrices, 
we fix these signs so that 
$\bar T=R(T)$, $\bar S_1=-R(S_1)$, $\bar S_2=-R(S_2)$, i.e.  
$$
\bar T=\left( 
\matrix 1 & 0 &  0 \cr
        1 & 1 &  0 \cr
        6 & 12 &  1 \cr \endmatrix \right),
\bar S_1=\left( 
\matrix 0 & 0 & -1 \cr
        0 & 1 &  0 \cr
       -1 & 0 &  0 \cr \endmatrix \right), 
\bar S_2=\left( 
\matrix -2 & -12 & -3 \cr
        1 & 5 &  1 \cr
       -3 & -12 &  -2 \cr \endmatrix \right).
$$ 
These matrices are understood as  
the monodromy around the regular singular points in the $x$-plane  
since the regular singularities of the 
Picard-Fuchs equation ({\PF}) are mapped to the corresponding 
points with non-trivial stabilizers under the map $t=t(x)$.  
The choices of the signs above are made so that they are consistent to 
the local properties of the PF equation, for example, 
$\text{det}(\bar S_1)=\text{det}(\bar S_2)=-1$ which come from 
the Picard-Lefshetz formula around the corresponding double points 
of $\tilde Q(\pm 6)$ and $\tilde Q(\pm 2)$ respectively.

\Proclaim{Proposition}{\Rgenerate} 
The matrices $R(T), R(S_1), R(S_2)$ generate the monodromy 
group $\Cal M_6(\vX)=O^+(U\oplus M_6)/\{ \pm \text{id}_{U\oplus M_6} \}$ 
defined in Definition {\sympMonod} $(n=6)$. 
\endproclaim
\demo{Proof} Note that the stabilizer $S_2$ may be related to 
the generators of $\Gamma_0(6)_+$ by 
$
S_2=\left( 
\matrix 
\sqrt{3} & {1\over \sqrt{3}} \cr
2\sqrt{3} & \sqrt{3} \cr 
\endmatrix \right) S_1 \;.
$
Then by Theorem {\Do}, 2) the claim follows. \qed 
\enddemo 

\Proclaim{Proposition}{\MS}  Consider $M_6=\langle 12 \rangle$, then one has: 
\item{1)} $\bar T,\bar S_1 \in O^+(U\oplus M_6)^*$. 
\item{2)} $\pm \bar S_2 \in O^+(U\oplus M_6) \setminus 
O^+(U\oplus M_6)^*$, i.e. 
$\pm \bar S_2$ does not come from the symplectic mapping class group 
of $\vX$. 
\item{3)} The monodromy representation of the symplectic diffeomorphisms 
$\Cal{MS}_6(\vX)(=O^+(U\oplus M_6)^*)$  is generated by 
$\bar T, \bar S_1, (\bar S_1 \bar S_2)^2$. 
\endproclaim
\demo{Proof} The discriminant is given by $A_{U\oplus M_6}\cong M_6^*/M_6 = 
{\bold Z \over 12} v / \bold Z v$. Since both $\bar T$ and 
$\bar S_1$ act on the basis $v$ by $v \mapsto v$ mod $U$, 
the induced action 
on the discriminant is the identity $\text{id}_{A_{U\oplus M_6}}$. 
This shows 1). 
On the other hand $\bar S_2$ acts on $v$ non-trivially, $v \mapsto 
5 v$ mod $U$. Therefore $\pm \bar S_2 \in O^+(U\oplus M_6)$ does not belong 
to the kernel $O^+(U\oplus M_6)^*$. By Theorem {\MainThII} we conclude 2). 

For 3), let us evaluate $(\bar S_1 \bar S_2)^2=(R(S_1)R(S_2))^2=
(R(S_2S_1))^2= R((S_2S_1)^2)$ ($R$ is an anti-homomorphism) and note that 
$$
S_2S_1=\left( 
\matrix 
\sqrt{3} & {1 \over \sqrt{3}} \cr
2\sqrt{3} & \sqrt{3} \cr 
\endmatrix \right)^2 
= \left( 
\matrix 
5 & 2 \cr
12 & 5 \cr 
\endmatrix
\right). 
$$
Therefore we obtain three generators $T, S_1, S_2S_1$ of $\Gamma_0(6)_{+6}$ 
which is isomorphic to $O^+(U\oplus M_6)^*$ by Theorem {\Do}, 1). 
\qed
\enddemo

\noindent
{\bf Remark.} The $\bar S_2$ represents the monodromy around $x={1\over 4}$ 
and has a simple geometric interpretation. 
Let us recall the explicit constructions of the algebraic cycles (lines) 
in $\tilde Q(\psi) (=\vX)$. As we have sketched in the proof of 
Theorem {\PeSt}, we may take 17 linearly independent lines in 
$NS(\tilde Q(\psi))$ which are independent of $\psi$ (coming from the 
toric divisors of the ambient space). The full lattice $NS(\tilde Q(\psi))$ 
is obtained by considering the lines, 
$m(\beta)$-lines, which are dependent 
of generic $\psi$ ($\beta+\beta^{-1}=\psi$). Explicitly we have 
$$
\beta_{\pm}={\psi \pm \sqrt{\psi^2-4} \over 2} 
= { 1 \pm \sqrt{1-4 x} \over 2 \sqrt{x} }  \qquad (x={1\over \psi^2}).
$$
Namely at $x={1\over 4}$ two lines $m(\beta_+)$ and $m(\beta_-)$ coincide, 
and these are exchanged under the monodromy operation around $x={1 \over 4}$. 
(The monodromy around $x=0$ is fictitious because of the isomorphism  
$\tilde Q(\psi)=\tilde Q(-\psi)$.) 
By Lemma {\LemmaI}, 1) the induced actions of the symplectic 
diffeomorphisms must be trivial on the lattice $NS(\vX)$ if the symplectic 
structure is generic.  Therefore the non-trivial monodromy behavior 
of the $m(\beta)$-lines explains that the monodromy $\bar S_2$ does 
not come from the symplectic diffeomorphism for a generic symplectic 
structure on $\vX$. We may argue similarly the monodromy 
around $x=\infty$.  \qed

\vskip0.5cm

Combining Proposition {\Rgenerate} and Proposition {\MS},  we evaluate   
directly the group index, $[\Cal M_n(\vX):\Cal{MS}_n(\vX)]=2$ $(n=6)$, 
which reproduce the number $p(6)=2$ of FM partners 
of the mirror $X$ (Proposition {\FMnumber}) from completely different 
picture.  Finally let us prove that the monodromy $\pm \bar S_2$ 
in Proposition {\MS}, 2) does not come from the autoequivalence of 
$D(X)$.

\Proclaim{Proposition}{\FMxyExpl} Let $X$ be a generic 
$M_6=\langle 12 \rangle$-polarizable K3 surface, i.e. 
the mirror of $\tilde Q(\psi)$, then 
the monodromy $\pm \bar S_2 \in O^+(U\oplus M_6)=O^+(\widetilde{NS}(X))$ 
does not come from $\text{Auteq} D(X)$.
\endproclaim
\demo{Proof} Since $\pm \bar S_2 \not\in O^+(U\oplus M_n)^*$, 
the element $(\pm \bar S_2,\text{id}_{U\oplus \vM_n}) 
\in O(U\oplus M_n)\times O(U\oplus \vM_n)$ does not extend to an 
isometry $O(\tilde \Lambda)$ ( where $U \oplus M_n \oplus U \oplus \vM_n 
\subset \tilde \Lambda$). Note that for $(X,\tau)$ generic, we have 
$\tau(T(X)) = U \oplus \vM_n$. Note also that when $\rho(X)=1$, 
the restriction $f\vert_{T(X)}$ of a Hodge isometry $f \in O_{Hodge}(\tilde H
(X,\bold Z), \bold C \omega_X)$ is trivial, i.e. $f\vert_{T(X)}=\pm 
\text{id}_{T(X)}$ (see for example [Og2,(4.1)]). 
Therefore $\bar S_2$ does not extend 
to a Hodge isometry in $O(\tilde \Lambda, \bold C \omega_1)$ 
$(\omega_1=\tau(\omega_X))$. By Theorem {\MainThI}, 
we obtain the desired result.  \qed
\enddemo

\noindent
{\bf Remark.} As explained in the previous Remark (after Proposition {\MS}), 
$\bar S_2 \in O^+(U\oplus M_6)$ comes from the monodromy about 
$x={1\over4}$, and the monodromy acts non-trivially on $NS(\vX)(=\vM_6)$ of 
$\vX=\tilde Q(\psi)$. Recall that the lines 
depicted in the toric diagram (the cube) in Fig.1 are independent of $\psi$, 
and generate a sublattice of rank 17 in $NS(\vX)$. This sublattice 
has signature $(1,16)$ ([PS,Proposition 1]). 
Therefore the monodromy around 
$x={1\over4}$ fixes this sublattice of signature $(1,16)$. From this 
argument, we see that the monodromy 
around $x={1\over4}$ preserves the orientations of the positive three 
planes in $U\oplus M_6 \oplus \vM_6 \subset \Lambda_{\text{K3}} \cong 
H^2(\vX,\bold Z)$. By the surjectivity 2) in Theorem 1.8, we may conclude 
that the monodromy does come from $Diff(\vX)$ although it does not come 
from $Symp(\vX,\kappa_\vX)$ (Proposition {\MS},2)).  
Now the mirror counterpart to this statement is clear. 
We expect that: {\it 
Let $X$ be a generic $M_6=\langle 12\rangle$-polarizable K3 surface and 
$Y \not\cong X$ be the FM partner of $X$. Then for the monodromy 
$\bar S_2$ there exists a kernel $\Cal E \in D(X \times Y)$ such that 
$ch(\Phi^{\Cal E}_{X \rightarrow Y})\vert_{U \oplus M_6}=\bar S_2$ in 
$O^+(U\oplus M_6)$. } 
Unfortunately, this does not follow directly from Theorem {\MainThI} 
(Main Theorem 1). 
\qed

\vfill\eject

\def\fsize#1{{\eightpoint\smc #1 }}
\def\efsize#1{{\eightpoint#1 }}
\def\spc{\hskip1cm}

\vskip1cm
\settabs 4 \columns   
\+ \spc \fsize{Shinobu Hosono$^\dagger$}         
                                    & & \fsize{Bong H. Lian}  \cr
\+ \spc\fsize{Department of mathematics}&&\fsize{Department of mathematics} \cr
\+ \spc \fsize{Harvard University}       && \fsize{Brandeis University}  \cr
\+ \spc \fsize{Cambridge, MA 02138}      && \fsize{Waltham, MA 02154}  \cr
\+ \spc \efsize{Email:}\efsize{hosono\@math.harvard.edu}       
                  &&  \efsize{Email:}\efsize{lian\@brandeis.edu} \cr
\+ \hskip2cm \efsize{hosono\@ms.u-tokyo.ac.jp}       &&  \fsize{}  \cr

\vskip1cm
\settabs 4 \columns   
\+ \spc \fsize{Keiji Oguiso$^\dagger$}         & & \fsize{Shing-Tung Yau}  \cr
\+ \spc\fsize{Department of mathematics}&&\fsize{Department of mathematics} \cr
\+ \spc \fsize{Harvard University}       && \fsize{Harvard University}  \cr
\+ \spc \fsize{Cambridge, MA 02138}      && \fsize{Cambridge, MA 02138}  \cr
\+ \spc \efsize{Email:}\efsize{keiji\@math.harvard.edu}       
                  &&  \efsize{Email:}\efsize{yau\@math.harvard.edu} \cr
\+ \hskip2cm \efsize{oguiso\@ms.u-tokyo.ac.jp}       &&  \fsize{}  \cr

\vskip0.3cm

{\leftskip0.5cm \eightpoint \noindent 
$^\dagger$ On leave of absence from Graduate School of Mathematical Sciences, 
University of Tokyo, \hfill\break 
Komaba 3-8-1, Meguro-ku, Tokyo 153-8914, JAPAN \par}

\enddocument